\documentclass[12pt]{article}

\usepackage{amsmath, amstext, amsgen, amsbsy, amsopn, amsfonts, amssymb, graphicx,pdfsync}

\newtheorem{theorem}{Theorem}

\newtheorem{definition}[theorem]{Definition}

\title{Links With Finite $n$-Quandles}
\author{Jim Hoste\\Pitzer College\\ \\ Patrick D. Shanahan\\Loyola Marymount University}
\begin{document}
\maketitle
\begin{abstract} We prove  a conjecture of Przytycki which asserts that the $n$-quandle of a link $L$ in the 3-sphere is finite if and only if the fundamental group of the $n$-fold cyclic branched cover of the 3-sphere, branched over $L$, is finite. 
\end{abstract}

\section{Introduction}\label{introduction}

While the algebraic study of racks and quandles dates back to the early 1900's, Fenn and Rourke in \cite{FR} credit Conway and Wraith with introducing the concepts in 1959 as an algebraic approach to study knots and links in 3-manifolds. In the late 1900's, several mathematicians began studying similar concepts under  names such as kei, distributive groupoids, crystals, and automorphic sets. In 1982, Joyce \cite{J} published a  ground-breaking work which included introducing the term quandle, giving both topological and algebraic descriptions of the fundamental quandle of a link, and proving that the fundamental quandle of a knot is a complete invariant  up to reversed mirror image. Much of Joyce's work was independently discovered by Matveev~\cite{M}.
In this article, we consider a quotient of the fundamental quandle of a link called the fundamental $n$-quandle, defined for any natural number $n$. Whereas the quandle of a link is usually infinite and somewhat untractable, there are many examples of knots and links for which the $n$-quandle is finite for some $n$.  In his Ph.D. thesis, Winker \cite{W} developed a method to produce the analog of the Cayley diagram for a quandle. In addition, Winker established a relationship between the $n$-quandle of the link $L$ and the fundamental group of $\widetilde{M_n}(L)$, the $n$-fold cyclic branched cover of the 3-sphere, branched over $L$.  When combined with previous work of Joyce, this implied that if the $n$-quandle of a link $L$ is finite, then so is $\pi_1(\widetilde{M_n}(L))$.  Przytycki \cite{P} then conjectured that this condition is both necessary and sufficient, which we prove to be true in this paper.  Our proof involves first generalizing a key result of Joyce: the cosets of the peripheral subgroup of a knot group can be given a quandle structure making it isomorphic to the fundamental quandle of the knot.  We extend this result to the $n$-quandle of a knot, showing that it can also be viewed as the set of cosets of the peripheral subgroup in a certain quotient of the knot group. This result  allows Winker's diagraming method to be replaced by the well known Todd-Coxeter method of coset enumeration.

We assume the reader is familiar with the theory of racks and quandles, but include basic definitions for completeness. The reader is referred to \cite{FR}, \cite{J},  \cite{J2}, \cite{M}, and \cite{W} for more information. A {\it  quandle} is a set $Q$ together with two binary operations $\rhd$ and $\rhd^{-1}$ which satisfies the following three axioms.
\begin{itemize}
\item[\bf Q1.] $x \rhd x =x$ for all $x \in Q$.
\item[\bf Q2.] $(x \rhd y) \rhd^{-1} y = x = (x \rhd^{-1} y) \rhd y$ for all $x, y \in Q$.
\item[\bf Q3.] $(x \rhd y) \rhd z = (x \rhd z) \rhd (y \rhd z)$ for all $x,y,z \in Q$.
\end{itemize}
A {\it rack} is more general, requiring only Q2 and Q3. It is important to note that, in general, the quandle operations are not associative. In fact, using axioms Q2 and Q3 it is easy to show that 
\begin{equation}\label{right distributive}
x \rhd ( y \rhd z) =  \left( (x \rhd^{-1} z) \rhd y \right) \rhd z.
\end{equation}
This property allows one to write any expression involving $\rhd$ and $\rhd^{-1}$ in a unique left-associated form (see \cite{W}). Henceforth, expressions without parenthesis are assumed to be left-associated. 

Given a quandle $Q$, each element $q\in Q$ defines a map $S_q:Q\to Q$ by $S_q(p)=p\rhd q$. It follows from axiom Q2 that $S_q$ is a bijection and $S_q^{-1}(p) = p \rhd^{-1} q$. From axiom Q3, it follows that $S_q$ is a quandle homomorphism. The automorphism $S_q$ is called the {\it point symmetry at $q$} and the set of all point symmetries generate the {\it inner automorphism group} $\text{Inn}(Q)$. A quandle $Q$ is {\it algebraically connected} if $\text{Inn}(Q)$ acts transitively on $Q$.  An {\it algebraic component} of $Q$ is a maximal algebraically connected subset of $Q$.

In \cite{J2}, Joyce defines two functors from the category of groups to the category of quandles.  These functors and their adjoints will be of importance in this paper. The first, denoted $\text{Conj}$, takes a group $G$ to a quandle $Q=\text{Conj}(G)$ defined as the set $G$ with operations given by conjugation. Specifically, $x\rhd y=y^{-1}xy$ and $x\rhd^{-1}y=yxy^{-1}$. Its adjoint, denoted $\text{Adconj}$ takes the quandle $Q$ to the group $\text{Adconj}(Q)$ generated by the elements of $Q$ and defined by the group presentation
$$\text{Adconj}(Q)=\langle \overline{q} \text{ for all $q$ in $Q$} \mid \overline{p\rhd q}=\overline{q}^{\,-1}\,\overline{p}\, \overline{q} \text{ for all $p$ and $q$ in $Q$}\rangle.$$

A quandle $Q$ is called an {\it $n$-quandle} if each point symmetry $S_q$ has order dividing $n$. It is convenient to write $x\rhd ^k y$ for $S^k_y(x)$, the $k$-th power of $S_y$ evaluated at $x$. Thus $Q$ is an $n$-quandle if for all $x$ and $y$ in $Q$, we have $x\rhd^ny=x$.
A second functor from groups to $n$-quandles  is defined for each natural number $n$ and is denoted $\text{Q}_n$.  Given a group $G$, the $n$-quandle $\text{Q}_n(G)$ is the set
$$\text{Q}_n(G)=\{x \in G \mid x^n=1\}$$
again with the operations given by conjugation. The adjoint of this functor is $\text{AdQ}_n$. If $Q$ is any $n$-quandle,  the group $\text{AdQ}_n(Q)$ is defined by the presentation
$$ \text{AdQ}_n(Q)=\langle \overline{q} \text{ for all $q$ in $Q$}  \mid \overline{q}^{\,n}=1, \overline{p\rhd q}=\overline{q}^{\,-1}\overline{p}\, \overline{q} \text{ for all $p$ and $q$ in $Q$}\rangle.$$

Quandles may be presented in terms of generators and relators in much the same way as groups. See \cite{FR} for a rigorous development of this topic. If the quandle $Q$  is given by the finite presentation
$$Q=\langle q_1, q_2, \dots, q_i\, | \, r_1, r_2, \dots ,r_j \rangle, $$
then Winker proves in \cite{W} that $\text{Adconj}(Q)$ and $\text{AdQ}_n(Q)$ can be finitely presented as
\begin{equation}\label{Winker presentation Adconj}
 \text{Adconj}(Q)=\langle \overline{q}_1, \overline{q}_2, \dots, \overline{q}_i \mid \overline{r}_1, \overline{r}_2, \dots, \overline{r}_j\rangle
\end{equation}
and
\begin{equation}\label{Winker presentation AdQn}
 \text{AdQ}_n(Q)=\langle \overline{q}_1, \overline{q}_2, \dots, \overline{q}_i \mid \overline{q}_1^{\,n}=1, \overline{q}_2^{\,n}=1, \dots, \overline{q}_i^{\, n}=1, \overline{r}_1, \overline{r}_2, \dots, \overline{r}_j\rangle.
 \end{equation}
Here, each quandle relation $r_i$ is an equation between two quandle elements each expressed using the generators, the operations $\rhd$ and $\rhd^{-1}$,  and parenthesis to indicate the order of operations. The associated group relation $\overline{r}_i$ must now be formed in a corresponding way  using conjugation. For example, if $r$ is the quandle relation $x=y\rhd(z\rhd^{-1}w)$, then $\overline{r}$ is the relation  $\overline{x}=\overline{w} \,\overline{z}^{\,-1}\,\overline{w}^{\, -1}\,\overline{y}\,\overline{w}\, \overline{z}\, \overline{w}^{\, -1}$.

\begin{figure}[h]
\vspace*{13pt}
\centerline{\includegraphics*[scale=.5]{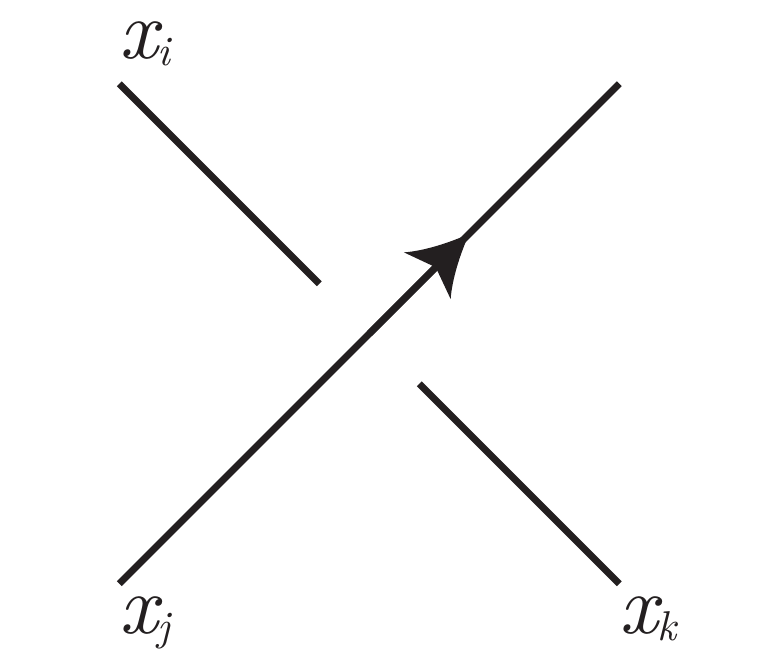}}
\caption{The relation $x_i=x_k \rhd x_j$ is associated to a  crossing with arcs labeled as shown.}
\label{crossing relation}
\end{figure}

Associated to every oriented knot or link $L$ in the 3-sphere $\mathbb S^3$ is its {\it fundamental quandle} $Q(L)$ which is defined by means of a presentation derived  from a regular diagram $D$ of $L$ with $a$ arcs and $c$ crossings. First assign quandle generators $x_1, x_2, \dots , x_a$ to each arc of $D$. Next, introduce a relation $r_\ell$ at each crossing of $D$ as shown in Figure~\ref{crossing relation}. It is easy to check that the three axioms, Q1, Q2, and Q3, are exactly what is needed to prove that $Q(L)$ is preserved by Reidemeister moves and hence is an invariant of the link. Passing from this presentation 
$$Q(L) = \langle x_1, \dots, x_a \, | \, r_1, \dots, r_c \rangle$$ 
to a presentation for $\text{Adconj}(Q(L))$  by using Winker's formula~\eqref{Winker presentation Adconj}, we obtain the well-known Wirtinger presentation of  $\pi_1(\mathbb S^3-L)$. Thus for any link $L$, $\pi_1(\mathbb S^3-L)\cong\text{Adconj}(Q(L))$.

Joyce proves in \cite{J} that $Q(L)$ is a complete invariant of  knots up to reverse mirror image. A less sensitive, but presumably more tractable, invariant is the {\it fundamental $n$-quandle} $Q_n(L)$ which can be defined for each natural number $n$.  If 
$$Q(L) = \langle x_1, \dots, x_a \, | \, r_1, \dots, r_c \rangle$$
is the presentation of the fundamental quandle of $L$ given by a diagram $D$ and $n$ is a fixed natural number, then the fundamental $n$-quandle of $L$ is defined to be the quandle with presentation
$$Q_n(L) = \langle x_1, \dots, x_a \, | \, r_1, \dots, r_c, s_1, \dots, s_k \rangle$$ 
where the relations $s_\ell$ are of the form $x_i \rhd^n x_j=x_i$ for all distinct pairs of generators $x_i$ and $x_j$.  As before, it is easy to check that $Q_n(L)$ is an invariant of $L$ and moreover that it is an $n$-quandle. Passing from this presentation of $Q_n(L)$ to a presentation for  $\text{AdQ}_n(Q_n(L))$  by using Winker's formula~\eqref{Winker presentation AdQn}, we see that $\text{AdQ}_n(Q_n(L))$ is a quotient of $\text{Adconj}(Q(L))$. In particular, we may present $\text{AdQ}_n(Q_n(L))$ by starting with the Wirtinger presentation of $\pi_1(\mathbb S^3-L)$ and then adjoining the relations $x^n=1$ for each Wirtinger generator $x$.  While the fundamental quandle of a nontrivial knot is always infinite, the associated $n$-quandle  is sometimes  finite. Determining when this occurs is the focus of this paper.

If $L$ is a link of more than one component, then both $Q(L)$ and $Q_n(L)$ are algebraically disconnected with one algebraic component $Q^i(L)$ and $Q^i_n(L)$, respectively, corresponding to each component $K_i$ of $L$. 

If $K$ is a knot, let $P$ be the peripheral subgroup of $G=\pi_1(\mathbb S^3-K)$ generated by the meridian $\mu$ and longitude $\lambda$ of $K$. In \cite{J}, Joyce defines a quandle structure on the set of right cosets $P\backslash G$ by declaring $Pg \rhd^{\pm 1} Ph=Pgh^{-1}\mu^{\pm 1} h$. He denotes this quandle as $(P\backslash G; \mu)$ and then proves that it is isomorphic to $Q(K)$.
This is the key step in Joyce's proof that the quandle is a complete knot invariant up to reverse mirror image.  It also implies that the order of $Q(K)$ is the index of $P$ in $G$ and hence that $Q(K)$ is infinite when $K$ is nontrivial. The key result of this paper is the following theorem which extends Joyce's result to the case of $Q_n(L)$. 
\begin{theorem}\label{key result} If $L=\{K_1, K_2, \dots, K_s\}$ is a link in $\mathbb S^3$ and $P_i$ is the subgroup of ${\rm AdQ}_n(Q_n(L))$ generated by the meridian $\mu_i$ and longitude $\lambda_i$  of $K_i$, then the quandle $(P_i\backslash {\rm AdQ}_n(Q_n(L)); \mu_i)$ is isomorphic to the algebraic component $Q^i_n(L)$ of $Q_n(L)$.  
\end{theorem}

Section 2 is devoted to proving Theorem~\ref{key result}. In Section 3 we use this result, as well as a theorem of Joyce, to prove the conjecture of Przytycki stated in the Abstract. Theorem~\ref{key result} implies that the Todd-Coxeter process for coset enumeration can be used to describe $Q^i_n(L)$ provided it is finite. In Section~4 we describe this in greater detail and give  examples.  In the last section, we enumerate all links that have finite $n$-quandles for some $n$. In a separate set of papers, we plan to  describe the $n$-quandles of these links, thereby providing a tabulation of all finite quandles that appear as the $n$-quandle of a link. The first of these papers is \cite{HS}, where we describe the 2-quandle of every Montesinos link of the form $M(p_1/2, p_2/2, p/q;e)$.  The authors extend their thanks to Daryl Cooper and Francis Bonahon for their assistance with Section~5. The authors also thank the referee for helpful comments.

\section{Relating {\boldmath $Q_n(L)$} to Cosets in {\boldmath $\text{AdQ}_n(Q_n(L))$}}

\begin{figure}[h]
\vspace*{13pt}
\centerline{\includegraphics*[scale=.5]{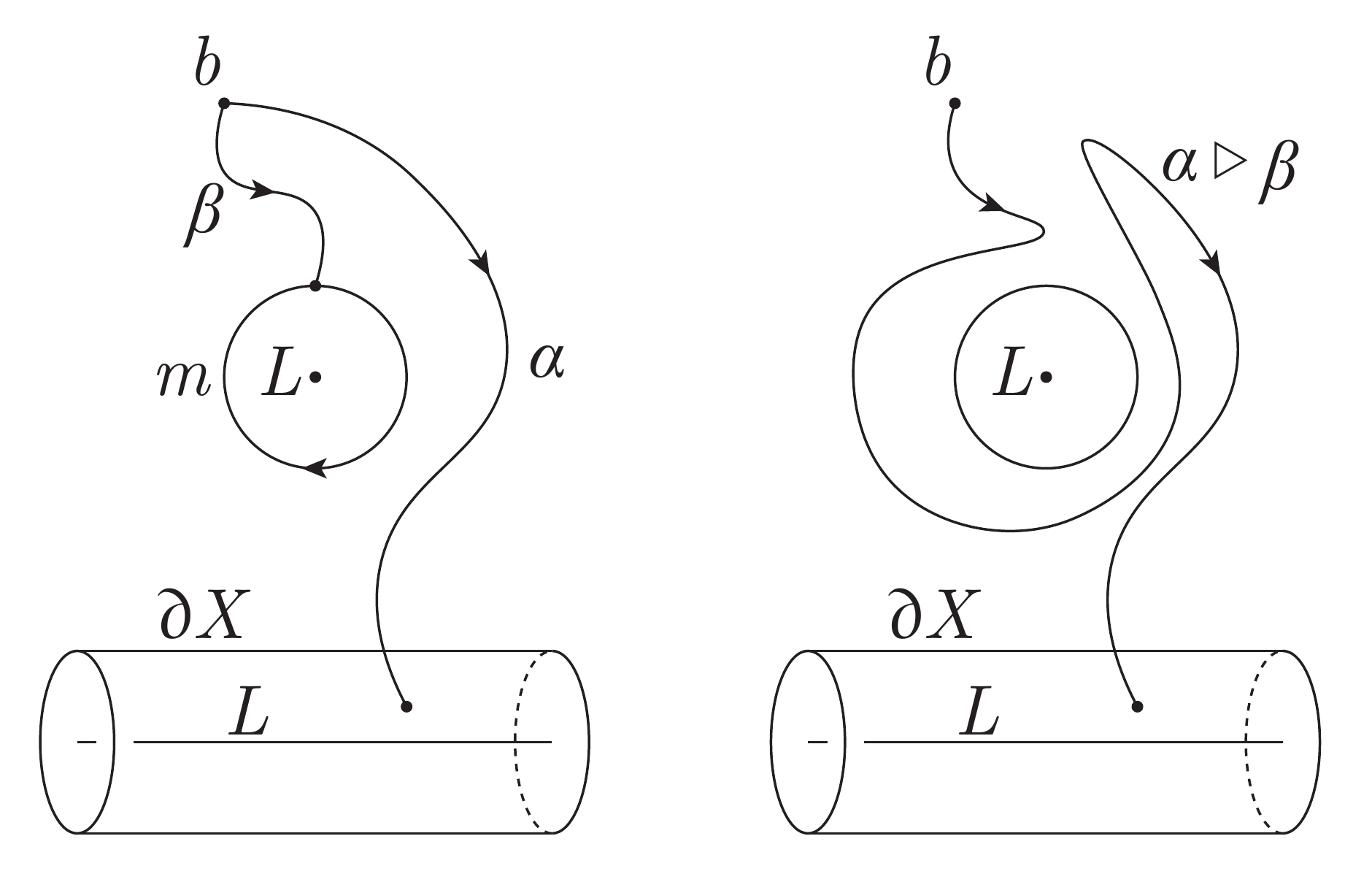}}
\caption{The topological definition of $\alpha \rhd \beta$.}
\label{topological defn}
\end{figure}

To prove Theorem~\ref{key result} we make use of  topological descriptions of both  the fundamental quandle  $Q(L)$ and  the $n$-quandle $Q_n(L)$.  We begin by recalling  Fenn and Rourke's formulation of $Q(L)$ given in \cite{FR}  and then extend it  to $Q_n(L)$. (Their formulation is actually for the rack associated to a framed link.) Let  $X=\mathbb S^3-\mathring{N}(L)$ be the exterior of $L$ and choose a basepoint $b$ in $X$. Define  $T(L)$ to be the set of all homotopy classes of  paths $\alpha:[0,1]\to X$ such that $\alpha(0)=b$ and $\alpha(1)\in \partial X$. Moreover,  we require that any homotopy be through a sequence of paths each of which starts at $b$ and ends at $\partial X$. Define the two binary operations, $\rhd$ and $\rhd^{-1}$, on  $T(L)$ by
\begin{equation}\label{topological quandle}
\alpha \rhd^{\pm 1} \beta=\beta m^{\mp 1}\beta^{-1}\alpha
\end{equation}
where $m$ is a meridian of $L$. Namely, $m$ is a loop in $\partial N(L)$ that begins and ends at $\beta(1)$, is essential in $\partial N(L)$,  is nullhomotopic in $N(L)$, and has linking number $+1$ with $L$. Thus the arc $\alpha \rhd \beta$ is formed by starting at the basepoint $b$, going along $\beta$ to $\partial N(L)$, traveling around $m^{-1}$, following $\beta^{-1}$ back to the base point, and  finally following $\alpha$ to its  endpoint in $\partial N(L)$.  See Figure~\ref{topological defn}. Note that the algebraic component  $T^i(L)$ corresponding to the $i$-th component $K_i$ of $L$ consists of those paths ending at $\partial N(K_i)$. The equivalence of $Q(L)$ and $T(L)$ is proven in  \cite{FR}. A similar description using ``nooses'' is given in \cite{J}.
In order to give a topological description of $Q_n(L)$ we introduce the following definition.
\begin{definition}\label{n equivalent topological}
\rm Suppose  $\alpha$  is a path in $X$  with $\alpha(0)=b$ and $\alpha(1) \in \{b\} \cup \partial X$. Suppose further that there exists $t_0$ with $0\le t_0\le 1$ such that $\alpha(t_0)\in \partial N(L)$. Let $\sigma_1(t)=\alpha(t t_0)$ and $\sigma_2(t)=\alpha((1-t)t_0+t)$. We say that  the path $\sigma_1 m^{\pm n}\sigma_2$ is obtained from $\alpha$ by a {\it $\pm n$-meridian move}. Two paths are called {\it $n$-meridionally equivalent} if they are related by a sequence of $\pm n$-meridian moves and homotopies.
\end{definition}

We  now define the $n$-quandle  $T_n(L)$ as the set of $n$-meridional equivalence classes of paths with the quandle operations defined by \eqref{topological quandle}.  Again, paths that end at $\partial N(K_i)$ give the algebraic component  $T^i_n(L)$ of $T_n(L)$.

\begin{theorem}\label{equivalence}
The $n$-quandles $Q_n(L)$ and $T_n(L)$ are quandle-isomorphic.
\end{theorem}

\noindent{\bf Proof.}  In \cite{FR}, the topological and algebraic-presentation definitions of the rack of a framed link are proven to be quandle isomorphic
by constructing homomorphisms $f:T \to Q$ and $g:Q \to T$ and then showing that both $f \circ g$ amd $g \circ f$ are the identity.  The same maps can be used to show that $T_n(L)$ and $Q_n(L)$ are isomorphic. Rather than repeating and extending Fenn and Rourke's proof here, we simply enumerate the differences from which the interested reader can easily fill in the details of the proof.
\begin{itemize}
\item In \cite{FR} homotopies in $T$ allow the endpoint of a path to move around on the chosen longitude of $L$ given by the framing, while we allow homotopies in $T_n$ to move the endpoint around in $\partial N(L)$. For our maps to be well-defined, this requires the idempotency axiom Q1 which is not present in a rack.
\item In $T_n$ we allow $n$-meridional moves that are not present in $T$.  In order for our maps to be well-defined this requires the addition of the corresponding relations $q_i\rhd^nq_j=q_i$ to $Q_n$. \hfill $\Box$
\end{itemize}

We are now prepared to prove Theorem~\ref{key result}.

{\renewcommand{\thetheorem}{\ref{key result}}
\begin{theorem} If $L=\{K_1, K_2, \dots, K_s\}$ is a link in $\mathbb S^3$ and $P_i$ is the subgroup of ${\rm AdQ}_n(Q_n(L))$ generated by the meridian $\mu_i$ and longitude $\lambda_i$  of $K_i$, then the quandle $(P_i\backslash {\rm AdQ}_n(Q_n(L)); \mu_i)$ is isomorphic to the algebraic component $Q^i_n(L)$ of $Q_n(L)$.  
\end{theorem}
\addtocounter{theorem}{-1}}

 \noindent{\bf Proof.}  Suppose that $L=\{K_1, K_2, \dots, K_s\}$. Without loss of generality, we shall prove the theorem for the first component $K_1$. 
 We begin by fixing some element $\nu \in Q_n(L)$ which we think of as a path from the basepoint $b$ in $X$ to   $\partial N(K_1)$. We now define  a map $\tau: \text{AdQ}_n(Q_n(L)) \to Q_n(L)$ by $\tau(\alpha)=\alpha^{-1}\nu$.
 
\noindent{\bf Claim 1:} The map $\tau$ is onto $Q^1_n(L)$.

{\bf Proof:} Let $\sigma$ be a path representing any element of $Q^1_n(L)$. Move $\sigma$ by a homotopy until $\sigma(1)=\nu(1)$ and let $\alpha$ be the loop $\alpha=\nu\sigma^{-1}$. Now $\tau(\alpha)=\alpha^{-1}\nu=\sigma\nu^{-1}\nu=\sigma$. 

Let $P_\nu$ be the subgroup of $\text{AdQ}_n(Q_n(L))$ generated by the meridian $\mu_1=\nu m \nu^{-1}$ and longitude $\lambda_1=\nu \ell \nu^{-1}$ of $K_1$. 

{\bf Claim 2:} $\tau^{-1}(\nu)=P_\nu$.

{\bf Proof:} Notice first that $\tau^{-1}(\nu)$ is a subgroup of  $\text{AdQ}_n(Q_n(L))$. For suppose that $\alpha, \beta \in \tau^{-1}(\nu)$. Now $\tau(\alpha\beta^{-1})=\beta \alpha^{-1} \nu=\beta \nu=\nu$ because $\alpha^{-1}\nu=\nu$ and $\beta^{-1}\nu=\nu$ implies $\nu=\beta\nu$.
Thus to show that $P_\nu \subset \tau^{-1}(\nu)$ we need only show that $\mu, \lambda \in \tau^{-1}(\nu)$. But $\tau(\lambda)=\lambda^{-1}\nu=(\nu \ell \nu^{-1})^{-1}\nu=\nu \ell^{-1}\nu^{-1} \nu=\nu \ell ^{-1}=\nu$ because $\ell \subset \partial X$. Similarly, $\mu\in \tau^{-1}(\nu)$. 

\begin{figure}[h]
\begin{tabular}{lll}
\includegraphics[scale=.4]{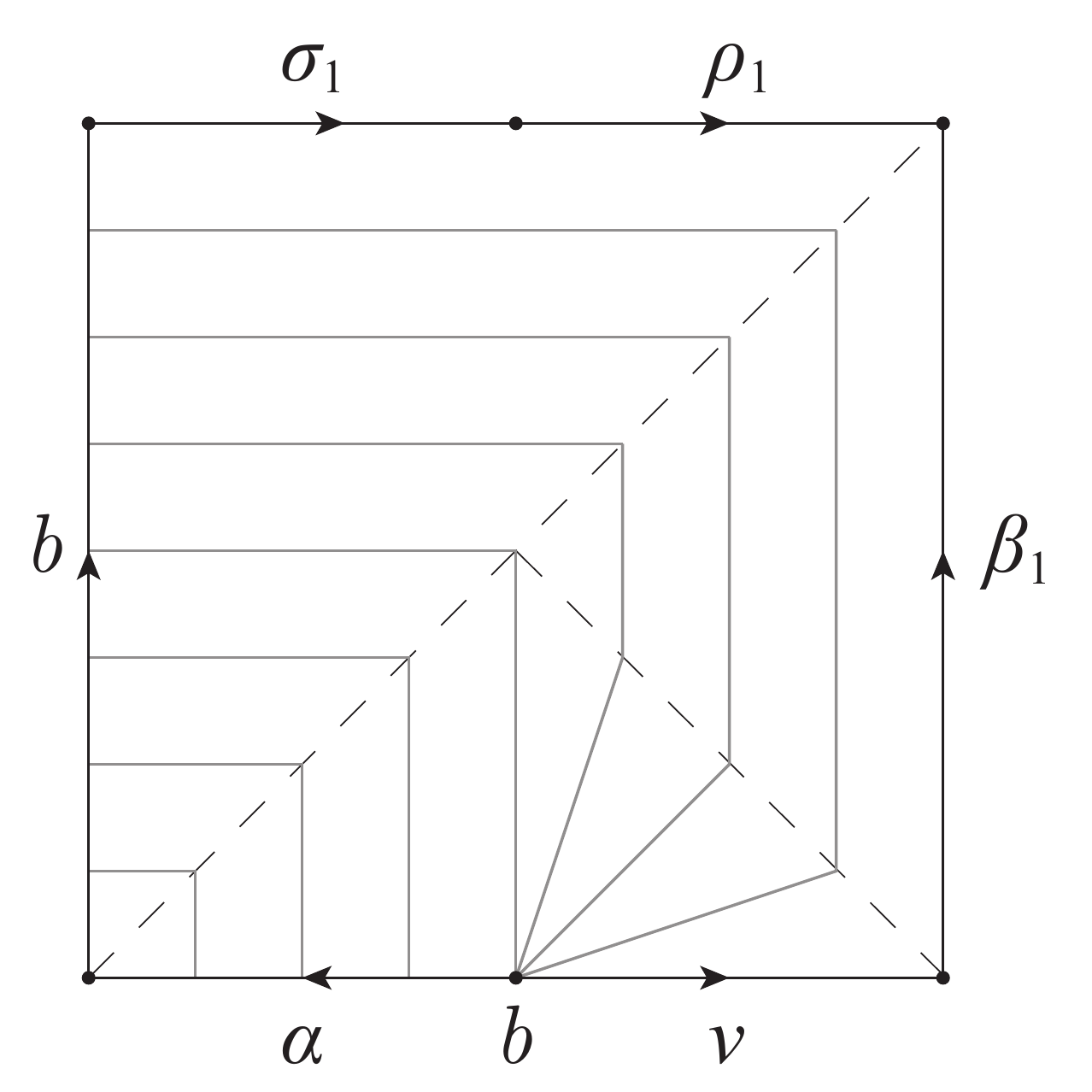}
&
\includegraphics[scale=0.4]{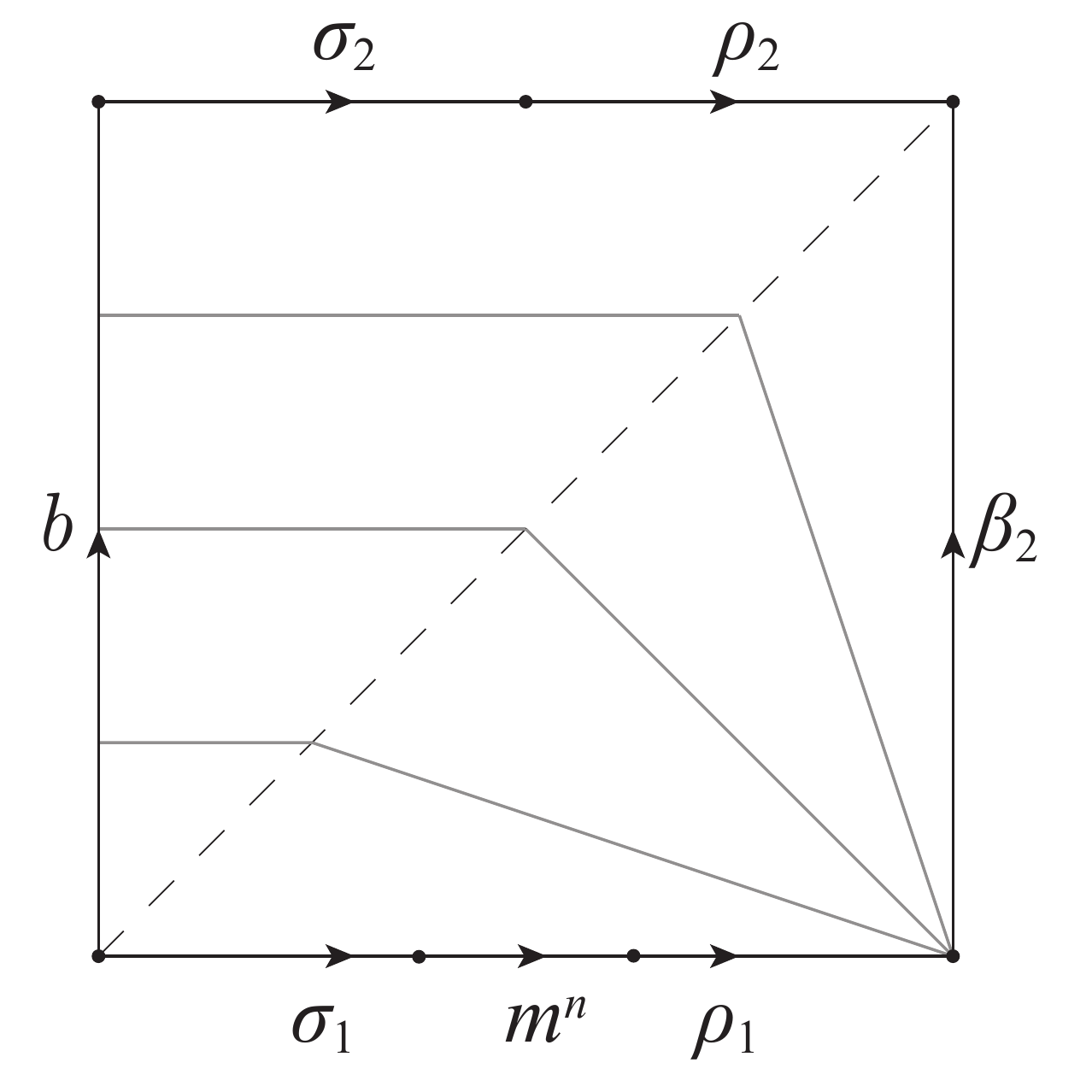}
&
\includegraphics[scale=0.4]{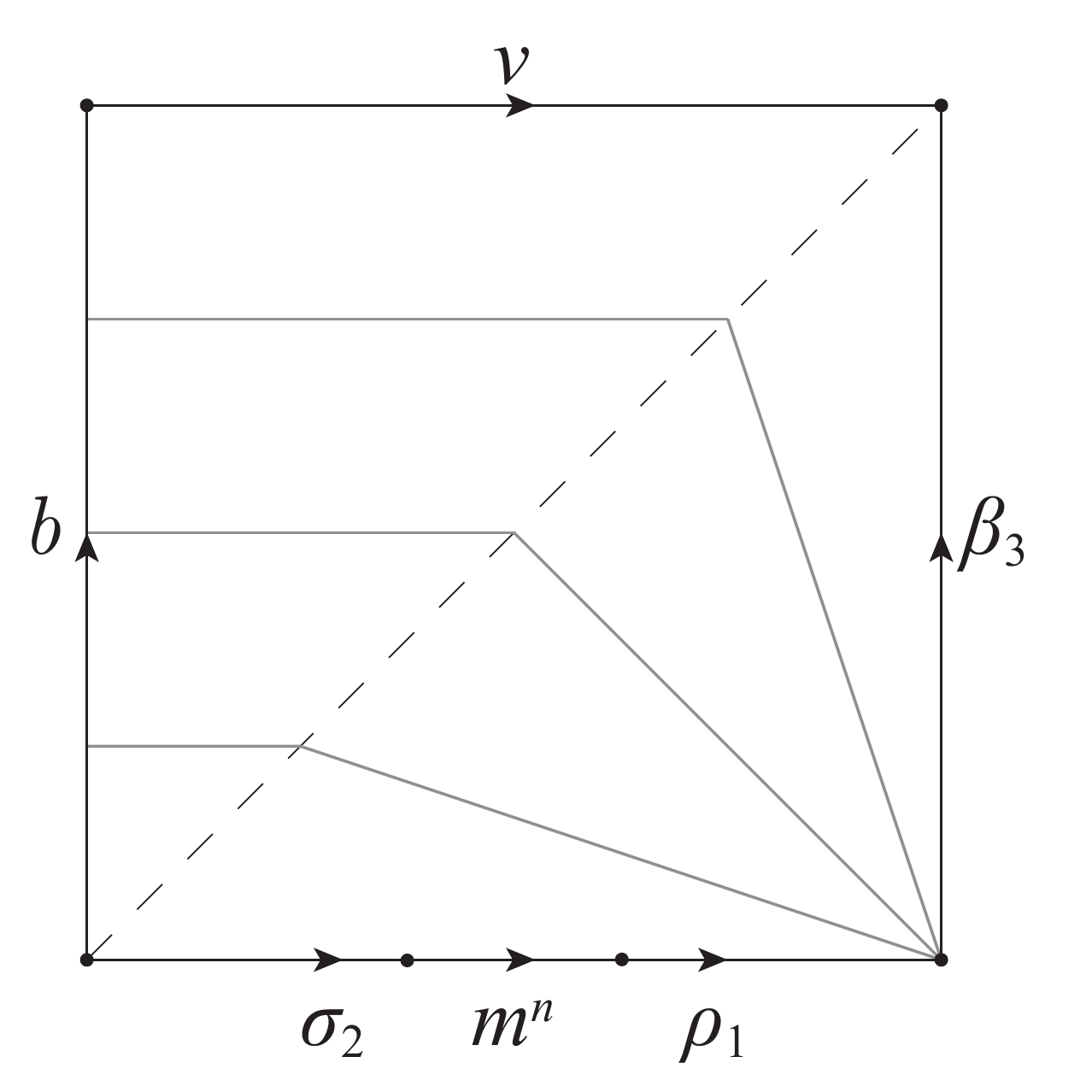}
\end{tabular}
\caption{Homotopies separated by $n$-meridian  moves. }
\label{homotopies}
\end{figure}

Now suppose that $\alpha \in \tau^{-1}(\nu)$. This means that $\alpha^{-1}\nu$ can be taken to $\nu$ by a sequence of $n$-meridian moves separated by homotopies. We illustrate the situation in Figure~\ref{homotopies}. The first homotopy begins at $\alpha^{-1}\nu$ and ends at the path $\sigma_1 \rho_1$ where $\sigma_1(1)=\rho_1(0)$ is a point in $\partial X$. We then do an $n$-meridian move, replacing $\sigma_1 \rho_1$ with the path $\sigma_1 m^{\pm n} \rho_1$. This path is then homotopic to the path $\sigma_2 \rho_2$ and so on until finally the last homotopy ends at $\nu$. For simplicity, the Figure illustrates the case of three homotopies separated by two $n$-meridian moves. Notice that the ``right edge'' of the $i$-th homotopy defines a path in $\partial N(K_1)$ which we call $\beta_i$. These homotopies can be reparameterized so that the polygonal paths indicated in each homotopy depict the new level sets. The first homotopy can now be thought of as one between the loop $\alpha$ and the loop $\nu \beta_1 \rho_1^{-1} \sigma_1^{-1}$. We then perform an $n$-meridian move to this loop and continue through the second homotopy, ending at the loop $\nu \beta_1 \beta_2 \rho_2^{-1} \sigma_2^{-1}$. Eventually we arrive at the loop $\nu \beta_1 \beta_2 \dots \beta_k \nu^{-1}$, an element of $P_\nu$. Thus $\alpha$  represents  an element of $P_\nu$ and hence $\tau^{-1}(\nu)\subset P_\nu$.

{\bf Claim 3:} Let $\phi_1$ be the automorphism of $\text{AdQ}_n(Q_n(L))$ given by conjugation by  $\mu_1$. Then  $\phi_1$ fixes every element of $P_\nu$.

{\bf Proof:} Suppose that $\nu\beta\nu^{-1}\in P_\nu$. Now
\begin{align*}
\phi_1(\nu\beta\nu^{-1})&=\mu_1^{-1}\nu\beta\nu^{-1}\mu_1\\
&=(\nu m\nu^{-1})^{-1}\nu\beta\nu^{-1}(\nu m\nu^{-1})\\
&=\nu m^{-1}\beta m \nu^{-1}\\
&=\nu m^{-1} m \beta \nu^{-1}\\
&=\nu \beta \nu^{-1}
\end{align*}
because loops in $\partial N(K_1)$ commute. 

We can now turn the set of right cosets $P_\nu \backslash \text{AdQ}_n(Q_n(L))$ into a quandle, which we denote as $(P_\nu \backslash \text{AdQ}_n(Q_n(L)); \mu_1)$ by defining
\begin{align}P_\nu \alpha \rhd^{\pm 1} P_\nu \beta&=P_\nu \phi_1^{\pm 1}(\alpha \beta^{-1})\beta\nonumber\\
&=P_\nu \mu_1^{\mp 1}\alpha \beta^{-1} \mu_1^{\pm 1} \beta \nonumber \\
&=P_\nu \alpha \beta^{-1} \mu_1^{\pm 1} \beta\label{quandle operation on cosets}
\end{align}
because $\mu_1 \in P_\nu$.

{\bf Claim 4:} The quandle operations defined in \eqref{quandle operation on cosets} are  well-defined.

{\bf Proof:} Suppose that $P_\nu \alpha=P_\nu a$ and $P_\nu \beta=P_\nu b$. Then
 \begin{align*}
\alpha \beta^{-1} \mu_1^{\pm 1}\beta (a b^{-1}\mu_1^{\pm 1} b)^{-1}&=\alpha \beta^{-1} \mu_1^{\pm 1}\beta  b^{-1}\mu_1^{\mp 1}ba^{-1}\\
&=\alpha \beta^{-1}\beta b^{-1}ba^{-1}\\
&=\alpha a^{-1} \in P_\nu
\end{align*}
because conjugation by $\mu_1^{\pm 1}$ fixes $\beta b^{-1}$, an element of $P_\nu$. Hence $P_\nu \alpha \rhd^{\pm1}P_\nu\beta=P_\nu a \rhd^{\pm1}P_\nu b$.

{\bf Claim 5:} The map $\tau$ determines a quandle isomorphism between  $(P_\nu \backslash \text{AdQ}_n(Q_n(L)); \mu_1)$ and  $Q^1_n(L)$.

{\bf Proof:} Define $\overline{\tau}:(P_\nu \backslash \text{AdQ}_n(Q_n(L)); \mu) \to Q^1_n(L)$ as $\overline{\tau}(P_\nu \alpha)=\tau(\alpha)$. Because $\tau^{-1}(\nu)=P_\nu$, it follows easily that $\overline{\tau}$ is both well-defined and injective. Because $\tau$ is onto $Q^1(L)$, we also have that $\overline{\tau}$ is onto  $Q^1_n(L)$. Thus $\overline{\tau}$ is a bijection. However, $\overline{\tau}$ is also a quandle homomorphism because 
\begin{align*}
\overline{\tau}(P_\nu \alpha \rhd P_\nu \beta)&= \overline{\tau}(P_\nu \alpha \beta^{-1}\mu_1^{-1}\beta)\\
&=\tau(\alpha \beta^{-1}\mu_1^{-1}\beta)\\
&=\beta^{-1}\mu_1^{-1}\beta\alpha^{-1}\nu\\
&=(\beta^{-1}\nu)m^{-1}(\beta^{-1}\nu)^{-1}(\alpha^{-1}\nu)\\
&=\tau(\alpha)\rhd \tau(\beta)\\
&=\overline{\tau}(P_\nu \alpha)\rhd \overline{\tau}(P_\nu \beta).
\end{align*}
\hfill $\square$

\pagebreak
\section{Przytycki's Conjecture}

In this section we prove the conjecture of Przytycki stated in the abstract.
\begin{theorem}\label{main theorem} Let $L$ be an oriented link  in  $\mathbb S^3$ and let $\widetilde{M_n}(L)$ be the $n$-fold cyclic branched cover of $\mathbb S^3$, branched over $L$. Then $Q_n(L)$  is finite,  if and only if  $\pi_1(\widetilde{M_n}(L))$ is finite.
\end{theorem}
Before giving the proof of Theorem~\ref{main theorem}, we point out the relationship between $\pi_1(\widetilde{M_n}(L))$ and a certain subgroup of $\text{AdQ}_n(Q_n(L))$. The reader is referred to \cite{W} for more details. If $M_n(L)$ is the $n$-fold cyclic cover of $\mathbb S^3-L$, then $\pi_1(M_n(L))$ is isomorphic to the subgroup $E^0$ of $\pi_1(\mathbb S^3-L)\cong \text{Adconj}(Q(L))$ consisting of those loops in $\mathbb S^3-L$ that lift to loops in the cover. Equivalently, $E^0$ consists of loops having total linking number zero with $L$, that is, those loops $\alpha$ such that the sum of the linking numbers of $\alpha$ with each component of $L$ is zero. The subgroup $E^0$ can also be described as those elements of  $\pi_1(\mathbb S^3-L)$ which, when written as words in the Wirtinger generators, have total exponent sum equal to zero. This concept is well-defined, and defines a subgroup, because each of the relators in the Wirtinger presentation has total exponent sum equal to zero. This last description extends to the quotient group $\text{AdQ}_n(Q_n(L))$.  Let $E^0_n$ be the subgroup of $\text{AdQ}_n(Q_n(L))$ consisting of all elements with total exponent sum equal to zero modulo $n$. In order to obtain the 
fundamental group of the cyclic branched cover we must algebraically kill the $n$-th power of each Wirtinger generator in $E^0$, hence,
\begin{equation}\label{Mn and En}
\pi_1(\widetilde{M_n}(L))\cong E^0_n.
\end{equation}
Notice further, that the index of $E^0_n$ in $\text{AdQ}_n(Q_n(L))$ is $n$.

One direction of Theorem~\ref{main theorem}  follows from work that appears in the  Ph.D. thesis of Joyce \cite{J2}. For completeness, and because  this result does not appear in Joyce's paper \cite{J}, we reproduce  his proof here (with some modification).

\begin{theorem}[Joyce]\label{Joyce's direction}  If $Q_n$ is any finite $n$-quandle, then $|{\rm AdQ}_n(Q_n)|\le n^{|Q_n|}$ and hence ${\rm AdQ}_n(Q_n)$ is finite.
\end{theorem}
\noindent{\bf Proof:} Suppose that $Q_n$ is a finite $n$-quandle with elements $\{q_1, q_2, \dots, q_k\}$. Now  $\text{AdQ}_n(Q_n)$ is generated by the ordered set of elements $\overline q_1, \overline q_2, \dots, \overline q_k$ so that every element in  $\text{AdQ}_n(Q_n)$ is a word in these generators and their inverses. 

{\bf Claim 1:} If $w=\overline q_{i_1}^{\,\epsilon_1}\overline q_{i_2}^{\,\epsilon_2}\dots\overline q_{i_m}^{\,\epsilon_m}$, where each exponent is $\pm 1$,  then we may rewrite $w$ as $w=\overline q_{j_1}^{\,\eta_1}\overline q_{j_2}^{\,\eta_2}\dots\overline q_{j_m}^{\,\eta_m}$, where each exponent is $\pm 1$, $j_1=\text{min}(j_1, j_2, \dots, j_m)$ and $j_1\le \text{min}(i_1, i_2, \dots, i_m)$.

\noindent{\bf Proof:} Suppose $\overline q_{i_k}^{\,\epsilon_k}$ is the first occurrence of the generator with smallest index and that $k>1$.  Now $q_{i_{k-1}}\rhd^{\epsilon_k} q_{i_k}=q_t$ for some $t$ and so $\overline q_{i_{k-1}}^{\,\epsilon_{k-1}}\overline q_{i_k}^{\,\epsilon_{k}}=\overline q_{i_k}^{\,\epsilon_{k}} \overline q_t^{\,\epsilon_{k-1}}$. If we replace  $\overline q_{i_{k-1}}^{\,\epsilon_{k-1}}\overline q_{i_k}^{\,\epsilon_{k}}$ with $\overline q_{i_k}^{\,\epsilon_{k}} \overline q_t^{\,\epsilon_{k-1}}$ in $w$, then either the first occurrence of the generator with smallest index has moved one place closer to the beginning of $w$, or a new generator of smaller index was introduced if $t<i_k$. Hence, after a finite number of steps of this kind, the first generator of $w$ will have the smallest index and it will be no greater than any of the indices in the original word.

{\bf Claim 2:} If $w=\overline q_{i_1}^{\,\epsilon_1}\overline q_{i_2}^{\,\epsilon_2}\dots\overline q_{i_m}^{\,\epsilon_m}$, where each exponent is $\pm 1$,  then we may rewrite $w$ as $w=\overline q_{j_1}^{\,\eta_1}\overline q_{j_2}^{\,\eta_2}\dots\overline q_{j_m}^{\,\eta_m}$, where each exponent is $\pm 1$ and $j_1\le j_2\le \dots \le j_m$.

\noindent{\bf Proof:} We proceed by induction on $m$. The case with $m=2$ is a direct consequence of Claim 1. Assume now that the result is true for words of length $m$ and suppose that 
$w=\overline q_{i_1}^{\,\epsilon_1}\overline q_{i_2}^{\,\epsilon_2}\dots\overline q_{i_{m+1}}^{\,\epsilon_{m+1}}$.  Applying the inductive hypothesis to the last $m$ generators of $w$, we may assume that $i_2\le i_3 \le \dots \le i_{m+1}$. If $i_1\le i_2$, we are done. If not, apply Claim 1 to $w$, which will strictly decrease the index of the first generator in $w$, and then  again apply the inductive hypothesis to the last $m$ generators. This cannot continue forever because the index of the first generator in $w$ cannot decrease below 1.

We may now write any word in $\text{AdQ}_n(Q_n)$ as $\overline q_{1}^{\,r_1}\overline q_{2}^{\,r_2}\dots\overline q_{k}^{\,r_k}$ and, using the fact that $\overline q_i^{\, n}=1$, we may assume that $0\le r_i<n$ for each $i$. There are at most $n^k=n^{|Q_n|}$
words of this kind. \hfill $\square$

\noindent {\bf Proof of Theorem~\ref{main theorem}:} Suppose $L$ is an oriented link and $Q_n(L)$ is finite. By Theorem~\ref{Joyce's direction}, it follows that 
$\text{AdQ}_n(Q_n(L))$ is finite. Hence  the subgroup $E^0_n$  of $\text{AdQ}_n(Q_n(L))$ is finite and so  
$\pi_1(\widetilde{M_n}(L))$ is finite  by \eqref{Mn and En}.

Now suppose that  $\pi_1(\widetilde{M_n}(L))$ is finite. Because $E^0_n$ has finite index in $\text{AdQ}_n(Q_n(L))$, it follows that $\text{AdQ}_n(Q_n(L))$ is finite.  Hence, for each component $K_i$ of $L$, the set of cosets $P_i\backslash\text{AdQ}_n(Q_n(L))$ is finite and therefore, by Theorem~\ref{key result},   each algebraic component $Q^i_n(L)$ of $Q_n(L)$ is finite. \hfill $\square$

\section{Examples} 

From the proof of Theorem~\ref{main theorem}, all information about the knot invariant $Q_n(L)$ is encoded  by the cosets of the subgroups $P_i$ in the group $\text{AdQ}_n(Q_n(L))$.  For example, if $Q_n(L)$ is finite, then
$$|Q_n(L)|=\sum_{i=1}^s \left[ \text{AdQ}_n(Q_n(L)) : P_i \right].$$
 Algorithmically computing the index of $P_i$ in the group $\text{AdQ}_n(Q_n(L))$ from a presentation of the group is a well-known problem in computational group theory. The first process to accomplish this task was introduced by Todd and Coxeter in 1936 \cite{TC} and  is now a  fundamental method in computational group theory. In addition to determining the index (if it is finite), the Todd-Coxeter process also provides a Cayley diagram that represents the action of right-multiplication on the cosets. In this section we will apply the Todd-Coxeter process to several examples and determine the quandle multiplication table from the Cayley diagram of the cosets.  More detailed treatments of the Todd-Coxeter process can be found in \cite{HO} and \cite{JO}.

Consider the right-hand trefoil knot $K$ and fix $n=3$.  From the Wirtinger presentation we obtain the presentation
$$ \text{AdQ}_3(Q_3(K))=\langle x, y \mid x^3=1, y^3=1, x^{-1}y^{-1}x y x y^{-1}=1\rangle.$$
A meridian for $K$ is $\mu =x$ and a (nonpreferred) longitude is $\lambda = yxxy$. The Todd-Coxeter process produces a {\em coset table} whose rows are numbered by indices $\alpha \in \{1,2,\dots, \kappa\}$ that represent cosets of $P$.  The columns are labeled by the generators and their inverses and encode the action of $\text{AdQ}_3(Q_3(K))$ on the cosets by right-multiplication.  An additional column will be added to give a representative $\phi(\alpha) \in \text{AdQ}_3(Q_3(K))$ of coset $\alpha$.

We initialize the coset table by letting $1$ represent the trivial coset $P$, thus $\phi(1)=e$ is a representative of this coset (we use $e$ here for the identity element of $\text{AdQ}_3(Q_3(K))$ to avoid confusion). Since $\mu = x \in P$, we have $Px=P$, this information is encoded in a {\em helper table} where $P$ is represented by index $1$ and is encoded in the coset table as a relation $1x=1$. Of course, it follows from this that $1x^{-1}=1$ as well, so there are two defined entries in row 1 of the coset table.

\centerline{
\raisebox{.05in}{\parbox{2in}{
$$\begin{array}{r|cccc|r}
 & x & y & x^{-1} & y^{-1} &  \phi \\ \hline 
1 & 1 &  & 1 &  & e  
\end{array}$$}}
\parbox{2in}{
$$\begin{array}{cc}
& x  \\ \hline
1 &  1
\end{array}$$}
}
Since $\lambda = yxxy \in P$ we also produce a helper table to encode $1yxxy=1$. Additional entries in the table are required to represent the cosets $1y$, $1yx$,  and $1yxx$. These entries are defined by adding indices $2$, $3$, and $4$, respectively, and adding additional information to the coset table for these indices coming from the helper table. For example, $2$ is defined to be the coset $1y$ and, thus, $1y=2$ and $2y^{-1}=1$ are encoded in the coset table.  At this point a {\em deduction} also occurs. Since $1yxxy=1$, we see in the helper table that $4y=1$.

\centerline{
\parbox{2in}{
$$\begin{array}{r|cccc|c}
 & x & y & x^{-1} & y^{-1} &  \phi \\ \hline 
1 & 1 & 2 & 1 & 4 & e  \\
2 & 3& & &1& y\\
3 & 4& & 2 & & yx\\
4 & &1 & 3 & & yx^2
\end{array}$$}
\parbox{2in}{
$$\begin{array}{ccccc}
& y & x & x & y \\ \hline
1 & 2 &3 &4 & 1
\end{array}$$}
}
This completes the initial set up of the coset table and is referred to as {\em scanning} the generators of $P$. The Todd-Coxeter process next proceeds to scan the relations of  $\text{AdQ}_3(Q_3(K))$ for all indices.
This encodes the fact that if $\alpha$ is any coset and $w=e \in \text{AdQ}_3(Q_3(K))$, then $\alpha w = \alpha$ in the coset table since $P \phi(\alpha) w = P \phi(\alpha)$ in $\text{AdQ}_3(Q_3(K))$. We scan the three relations $x^3=e$, $y^3=e$, and $x^{-1} y^{-1} x y x y^{-1}=e$, in this order, for each index, defining new indices and obtaining new deductions along the way.

Scanning $x^3$ for $\alpha =1$ gives no new information. Scanning $y^3$ gives no new definitions but does produce the deduction $2y=4$ and scanning $x^{-1} y^{-1} x y x y^{-1}$ defines the indices 5 and 6 as shown in the coset tables below.

 \centerline{
\parbox{2in}{
$$\begin{array}{r|cccc|c}
 & x & y & x^{-1} & y^{-1} &  \phi \\ \hline 
1 & 1 & 2 & 1 & 4 & e  \\
2 & 3& 4& &1& y\\
3 & 4& & 2 & & yx\\
4 & &1 & 3 &2 & yx^2
\end{array}$$}
\parbox{2in}{
$$\begin{array}{cccc}
& y & y & y \\ \hline
1 & 2 &4 &1 
\end{array}$$}
}

\centerline{
\parbox{2in}{
$$\begin{array}{r|cccc|c}
 & x & y & x^{-1} & y^{-1} &  \phi \\ \hline 
1 & 1 & 2 & 1 & 4 & e  \\
2 & 3& 4& 6&1& y\\
3 & 4& & 2 & & yx\\
4 & 5&1 & 3 &2 & yx^2 \\
5 & & 6&4 & & yx^3\\
6 &2 & & & 5 & yx^3y
\end{array}$$}
\hspace{.2in}
\parbox{2in}{
$$\begin{array}{ccccccc}
& x^{-1} & y^{-1} & x & y & x & y^{-1} \\ \hline
1 &1 &4 &5 &6 &2 & 1
\end{array}$$}
}
At this point we see that the representative for coset 5 is $\phi(5)=yx^3$. Since $x^3=e$ in the group $\phi(5)=yx^3=y=\phi(2)$ and so the cosets $5$ and $2$ are the same. This information is determined by a {\em coincidence} which occurs when scanning $x^3$ for $\alpha = 2$.
Filling in the entries of the helper table from left to right, $2x=3$, $3x=4$, $4x=5$. However we require $2xxx=2$ thus we see that $5=2$. In the coset table we process this coincidence by replacing all values of $5$ with $2$, merging the data from row $5$ into row $2$, and then deleting row $5$. In merging the data from $5$ to $2$ we see a new coincidence, namely $6=4$ and so we repeat the coincidence procedure for $6=4$ before moving on to the next scan.

\centerline{
\parbox{2in}{
$$\begin{array}{r|cccc|c}
 & x & y & x^{-1} & y^{-1} &  \phi \\ \hline 
1 & 1 & 2 & 1 & 4 & e  \\
2 & 3& 4& \not \!6 4&1& y\\
3 & 4& & 2 & & yx\\
4 & \not \! 5 2&1 & 3 &2 & yx^2 \\
\not \!5 & & \not \!6&\not \!4 & & yx^3\\
\not \!6 &\not \!2 & & & \not \! 5\!\not\!2 & yx^3y
\end{array}$$}
\parbox{2in}{
$$\begin{array}{cccc}
& x & x & x \\ \hline
2 & 3 & 4 & 5=2
\end{array}$$}
}
Scanning $x^{-1} y^{-1} x y x y^{-1}$ for $\alpha=2$ completes the table. The process terminates  after the table is complete and all relations have been scanned for all indices. In our example, no additional coincidences occur and the completed table is shown below.
\begin{table}[h]
$$\begin{array}{r|cccc|c}
 & x & y & x^{-1} & y^{-1} &  \phi \\ \hline 
1 & 1 & 2 & 1 & 4 & e  \\
2 & 3& 4& 4&1& y\\
3 & 4& 3& 2 &3 & yx\\
4 & 2&1 & 3 &2 & yx^2
\end{array}$$
\caption{Completed coset table for $P\backslash\text{AdQ}_3(Q_3(K))$ where $K$ is the trefoil knot}
\end{table}

It is important to note that the operation encoded by the coset table is that of right-multiplication. It is not the operations of $\rhd^{\pm 1}$ in the quandle $P\backslash\text{AdQ}_3(Q_3(K))$.  The multiplication table for the quandle can be easily worked out, however, from the coset table and the definition of the operations $P g \rhd^{\pm 1} P h = P g h^{-1} x^{\pm 1} h$ since $\mu =x$.  From the completed coset table, the quandle $Q_3(K)$ has four elements $P$, $Py$, $Pyx$, and $Pyx^2$.  So, for example, $Py \rhd Pyx = Py x^{-1} y^{-1} x y x$. This coset is represented by $1y x^{-1} y^{-1} x y x = 4$ in the coset table. Therefore, $Py \rhd Pyx = Pyx^2$.  The full multiplication table for $Q_3(K)$ is given below.

\begin{table}[h]
$$\begin{array}{l|cccc}
\rhd & P & Py & Pyx & Pyx^2 \\ \hline 
P & P & Pyx^2& Py & Pyx\\
Py & Pyx & Py & Pyx^2& P \\
Pyx & Pyx^2 &P & Pyx & Py \\
Pyx^2 & Py &Pyx & P& Pyx^2 
\end{array}$$
\caption{The multiplication table for $Q_3(K)$ where $K$ is the trefoil knot.}
\end{table}

Applying the Todd-Coxeter method in the case of the trefoil for $n=2, 3, 4, 5$, enumerating the cosets of both the trivial subgroup as well as $P=\langle \mu, \lambda \rangle$, we obtain the  data in Table~\ref{trefoil data}. These calculations agree with the  well known fact that  $\pi_1(\widetilde{M_n}(K))$ for the trefoil with $n=2, 3, 4$, or $5$ is, respectively, the cyclic group of order 3, the quaternion group of order 8, the binary tetrahedral group of order 24, and the binary icosahedral group of order 120. See \cite{R}.

\begin{table}[h]
\begin{center}
\begin{tabular}{ccccc}
$n$&$|P|$&$|Q_n(K)|$&$|\text{AdQ}_n(Q_n(K))|$& $|\pi_1(\widetilde{M_n}(K))|$\\
\hline
2&2&3&6&3\\
3&6&4&24&8\\
4&16&6&96&24\\
5&50&12&600&120
\end{tabular}
\end{center}
\label{trefoil data}
\caption{The order of $\text{AdQ}_n(Q_n(K))$ and index of $P$ for the right-handed trefoil.}
\end{table}

As another example, consider the $(2, 2, 3)$-pretzel link $L$ and fix $n=2$. Starting with the standard pretzel diagram with Wirtinger generators $x, y, z$, we obtain the following presentation of $\text{AdQ}_2(Q_2(L))$.
\begin{align*} \langle x, y, z \, \mid \, & x^{-1}z^{-1}xzxy^{-1}x^{-1}y=1, y^{-1}x^{-1}yx y z y z^{-1}y^{-1}z^{-1}=1,\\
& y^{-1}x^{-1}yxyzyz^{-1}y^{-1}x^{-1}y^{-1}xyx^{-1}z^{-1}z=1, x^2=1, y^2=1, z^2=1 \rangle
\end{align*}
The link $L$ has two components and subgroups generated by a meridian and longitude of each component are $P_1=\langle x, x^{-1}zxy^{-1}\rangle$ and $P_2=\langle y, y^{-1}x^{-1}yzyz^{-1}xzyzy^{-1} \rangle$.  Applying the Todd-Coxeter process for each of these subgroups gives 
$$\left| Q_2(L) \right| = \left[ \text{AdQ}_2(Q_2(L)) : P_1 \right] + \left[ \text{AdQ}_2(Q_2(L)) : P_2 \right] = 8+24 = 32.$$
These calculations agree with Theorem~1.1 of \cite{HS} where it is shown using Winker's diagramming method \cite{W} that if $L$ is the Montesinos link of the form $(1/2, 1/2, p/q; e)$, then $|Q_2(L)|=2(q+1)|(e-1)q-p|$. For the $(2, 2, 3)$-pretzel link $L$ we have $p=1$, $q=3$, and $e=0$.

\section{Links With Finite $n$-Quandles}

The set of links which have a finite $n$-quandle for some $n$ can be derived from Thurston's geometrization theorem.  To see this, let $L$ be a link and $n>1$ an integer such that $Q_n(L)$ is finite. By Theorem~\ref{main theorem}, we have that $\pi_1(\widetilde{M_n}(L))$ is finite. Define  $\mathcal O (L,n)$ to be the 3-orbifold with underlying space $\mathbb S^3$ and singular locus $L$ where each component of $L$ is labelled $n$.  (Both \cite{BP} and \cite{CH} are excellent references for orbifolds.)  We now have  a manifold covering of the orbifold, $p: \widetilde{M_n}(L) \rightarrow \mathcal O (L,n)$, and the covering map $p$ induces a homomorphism $p_* : \pi_1(\widetilde{M_n}(L)) \rightarrow \pi_1^{orb}(\mathcal O (L,n))$ for which the index of $p_*(\pi_1(\widetilde{M_n}(L)))$ in $\pi_1^{orb}(\mathcal O (L,n))$ is the branch index $n$. Since $\pi_1(\widetilde{M_n}(L))$ is finite, it follows that $\pi_1^{orb}(\mathcal O (L,n))$ is finite. In addition, the universal orbifold cover of $\mathcal O(L,n)$ is a simply-connected manifold (equal to the universal cover of $\widetilde{M_n}(L)$) and, since $\pi_1^{orb}(\mathcal O (L,n))$ is finite, the universal cover is also compact. Now Thurston's geometrization theorem asserts that the only compact, simply-connected 3-manifold is $\mathbb S^3$.  Therefore, $\mathcal O(L,n)$ is a spherical 3-orbifold.  In \cite{DU}, Dunbar classifies all geometric, non-hyperbolic 3-orbifolds.  The following, obtained from Dunbar, is the complete list of all spherical 3-orbifolds with underlying space $\mathbb S^3$ and singular locus $L$ with each component labelled $n$. Therefore, it also represents the list of all links in $\mathbb S^3$ with finite $Q_n(L)$ for some $n$.

In Table~\ref{list of links}, we list the links as they appear in \cite{DU}. A box labeled $k$ denotes $k$ left-handed half twists between the two strands and a box labeled $m/n$ denotes the $m/n$ rational tangle with $-n/2 \le m \le n/2$ and $m\ne 0$. See \cite{DU} for a detailed explanation. 

\begin{table}
\label{list of links}
$$
\begin{array}{ccccc}
\includegraphics[width=1.25in,trim=0 50pt 0 75pt,clip]{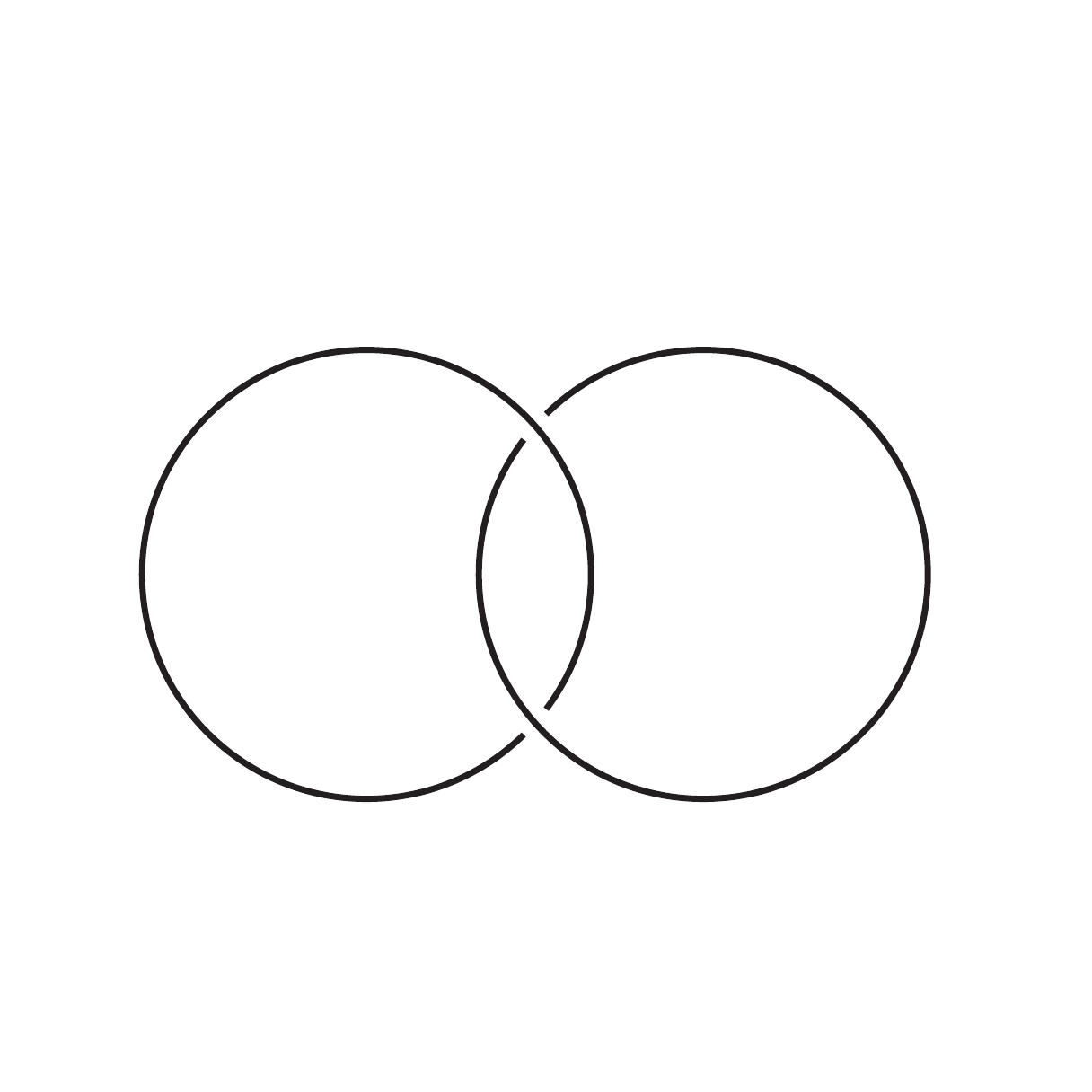} &  & \includegraphics[width=1.5in,trim=0 60pt 0 65pt,clip]{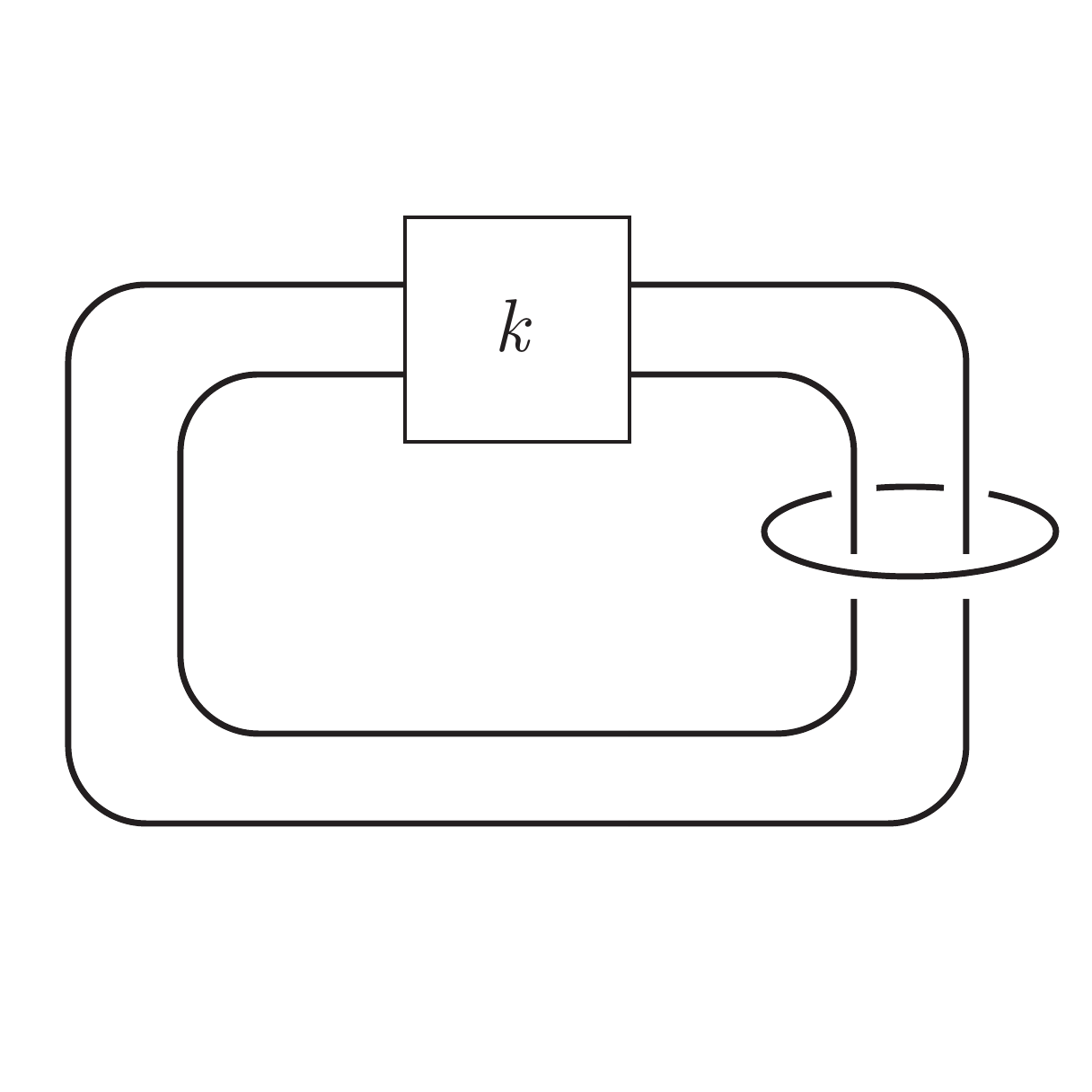} & & \includegraphics[width=1.25in,trim=0 50pt 0 50pt,clip]{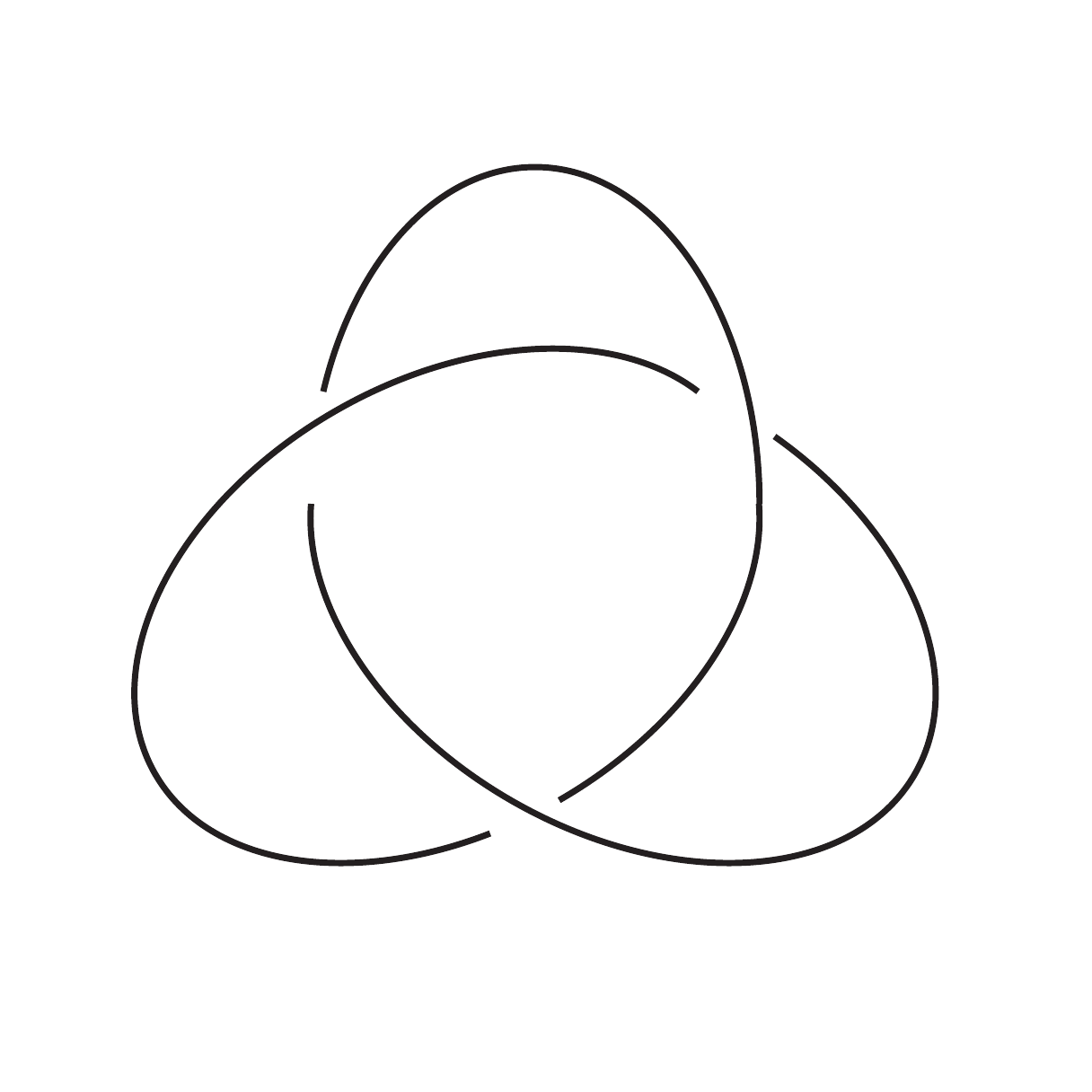} \\
\scriptstyle n > 1 & & \scriptstyle k \neq 0,\  n=2& & \scriptstyle n=3, 4, 5 \\
\\
\includegraphics[width=1.25in,trim=0 25pt 0 0,clip]{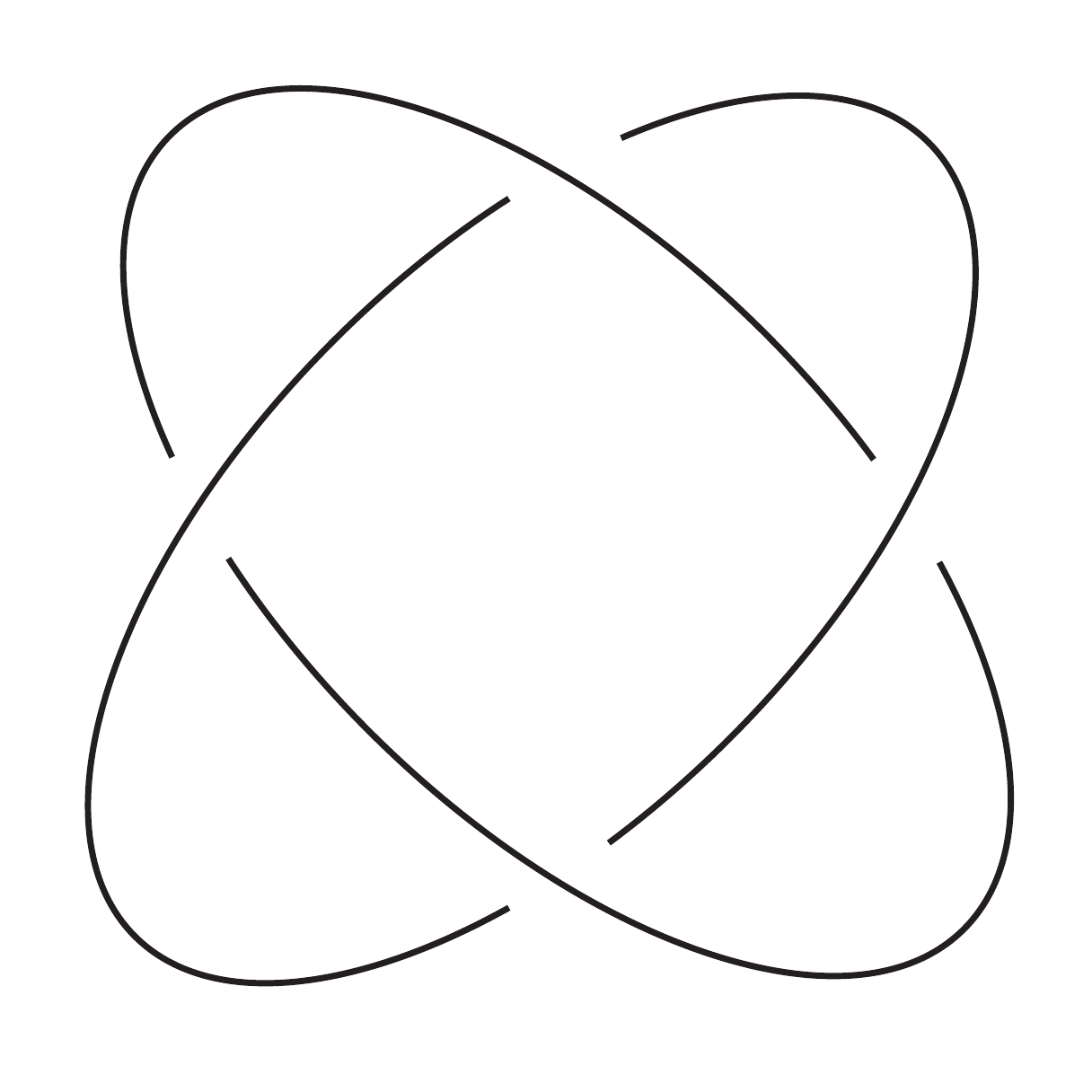} & & \includegraphics[width=1.25in,trim=0 25pt 0 0,clip]{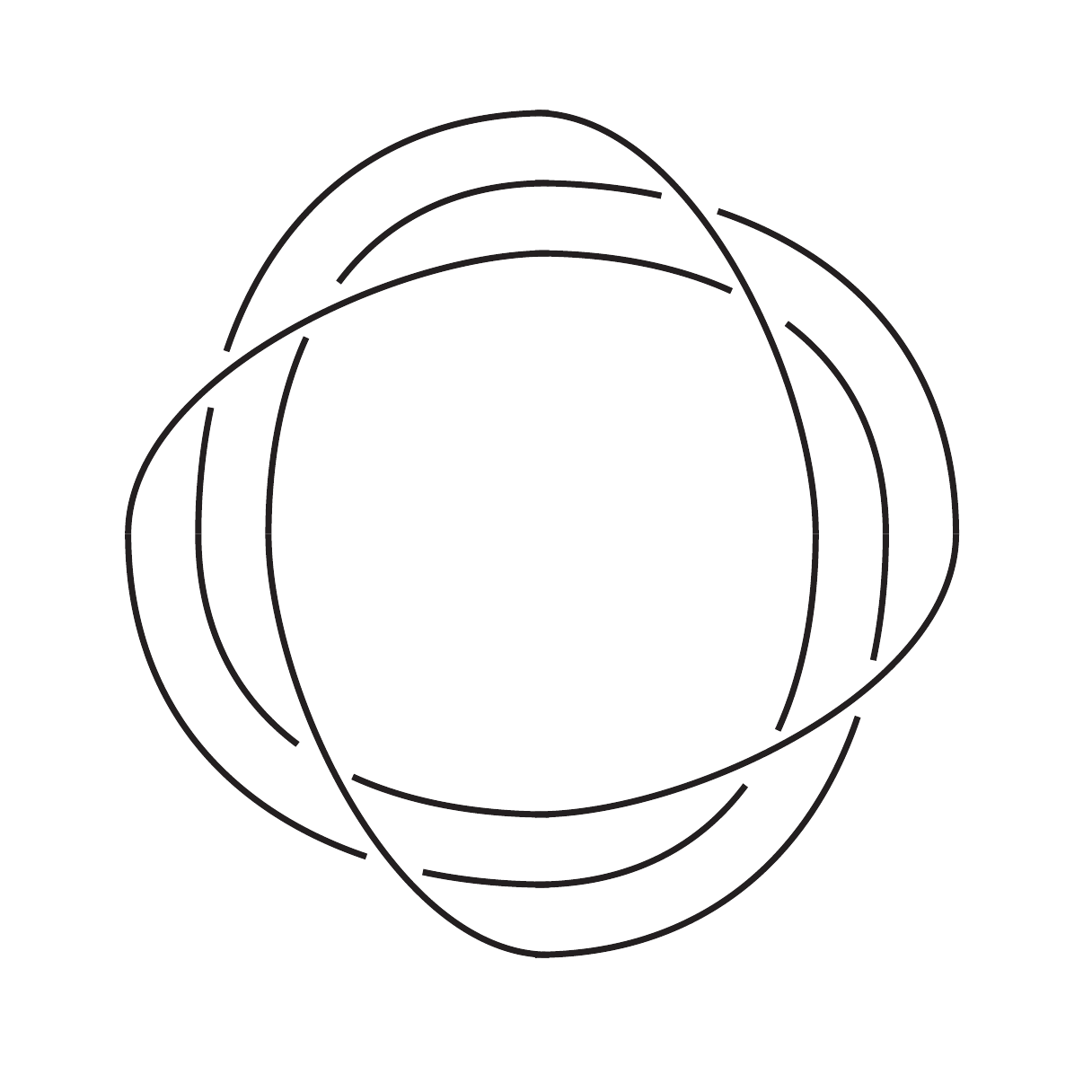} & & \includegraphics[width=1.25in,trim=0 25pt 0 0,clip]{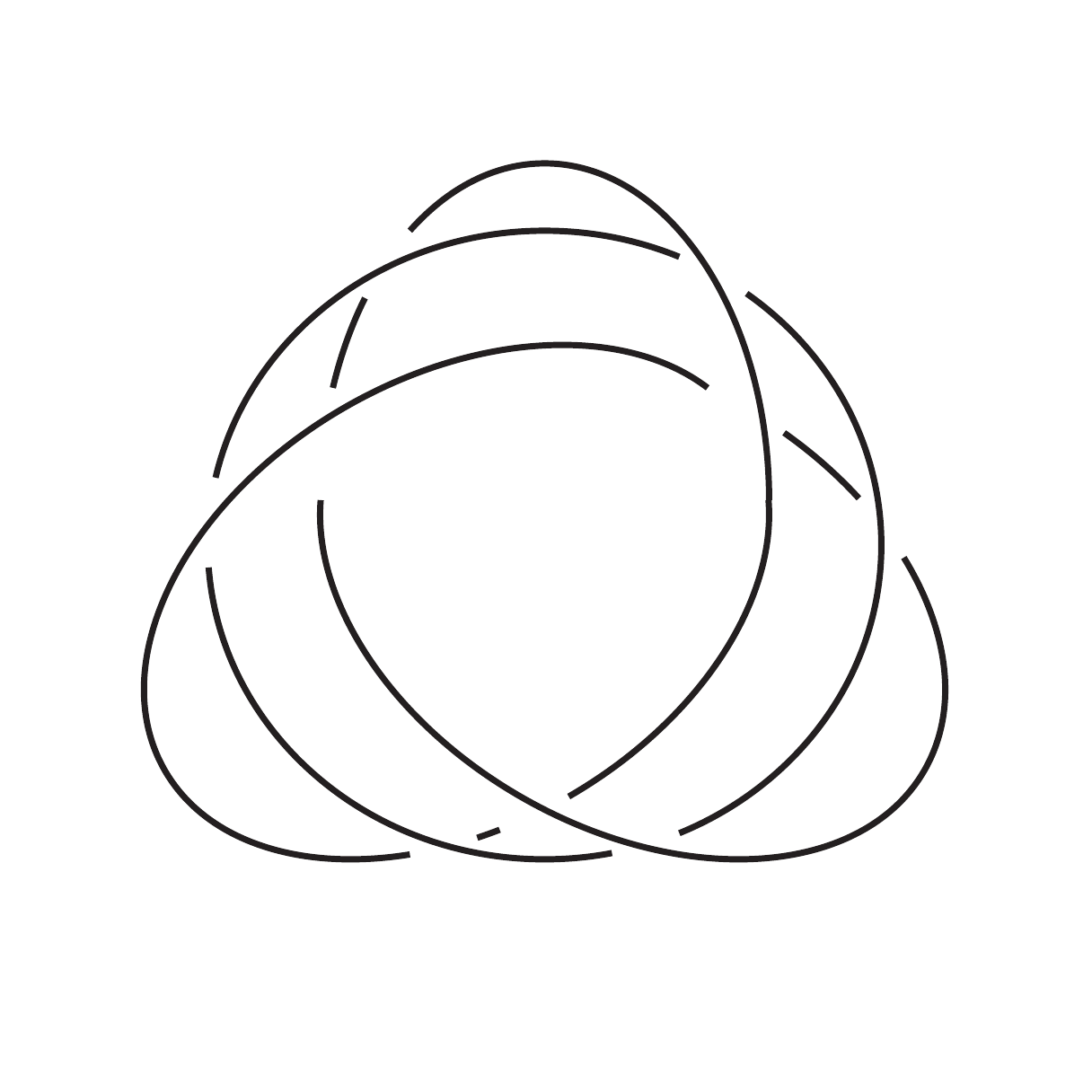} \\
\scriptstyle n =3 & &\scriptstyle n=2& &\scriptstyle n=2 \\
\\
\includegraphics[width=1.25in,trim=0 25pt 0 0,clip]{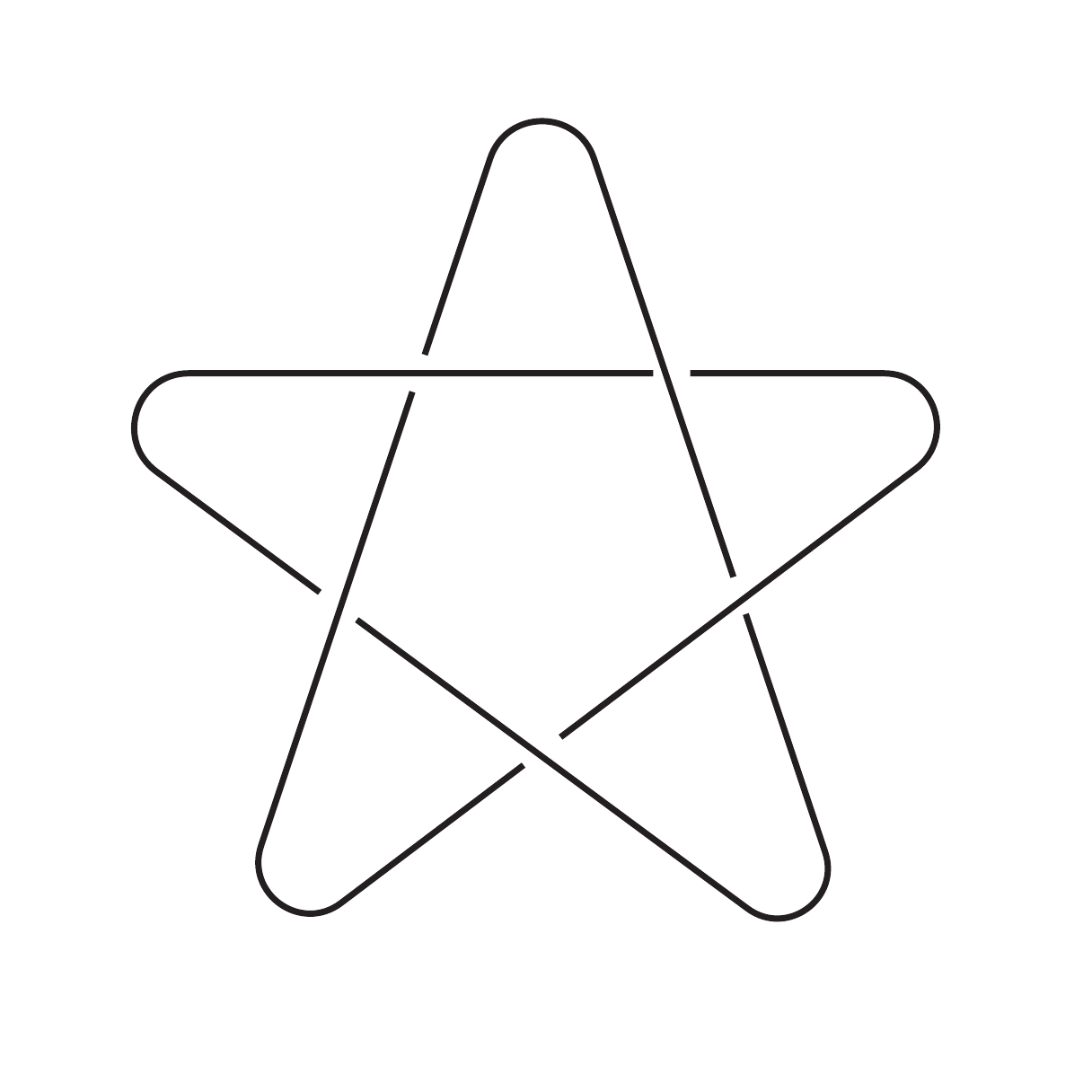} & & \includegraphics[width=1.25in,trim=0 25pt 0 0,clip]{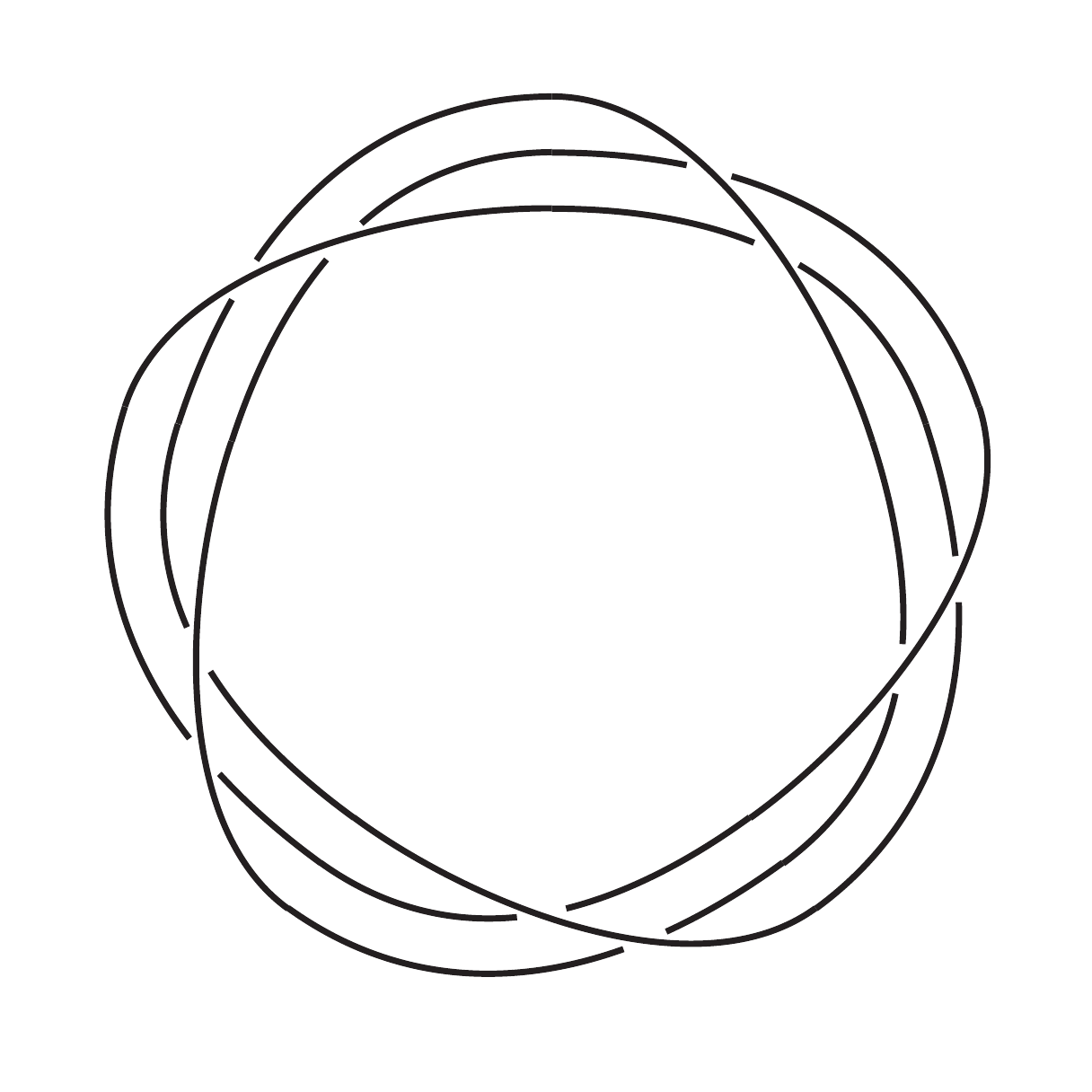} & & \includegraphics[width=1.5in,trim=0 50pt 0 25pt,clip]{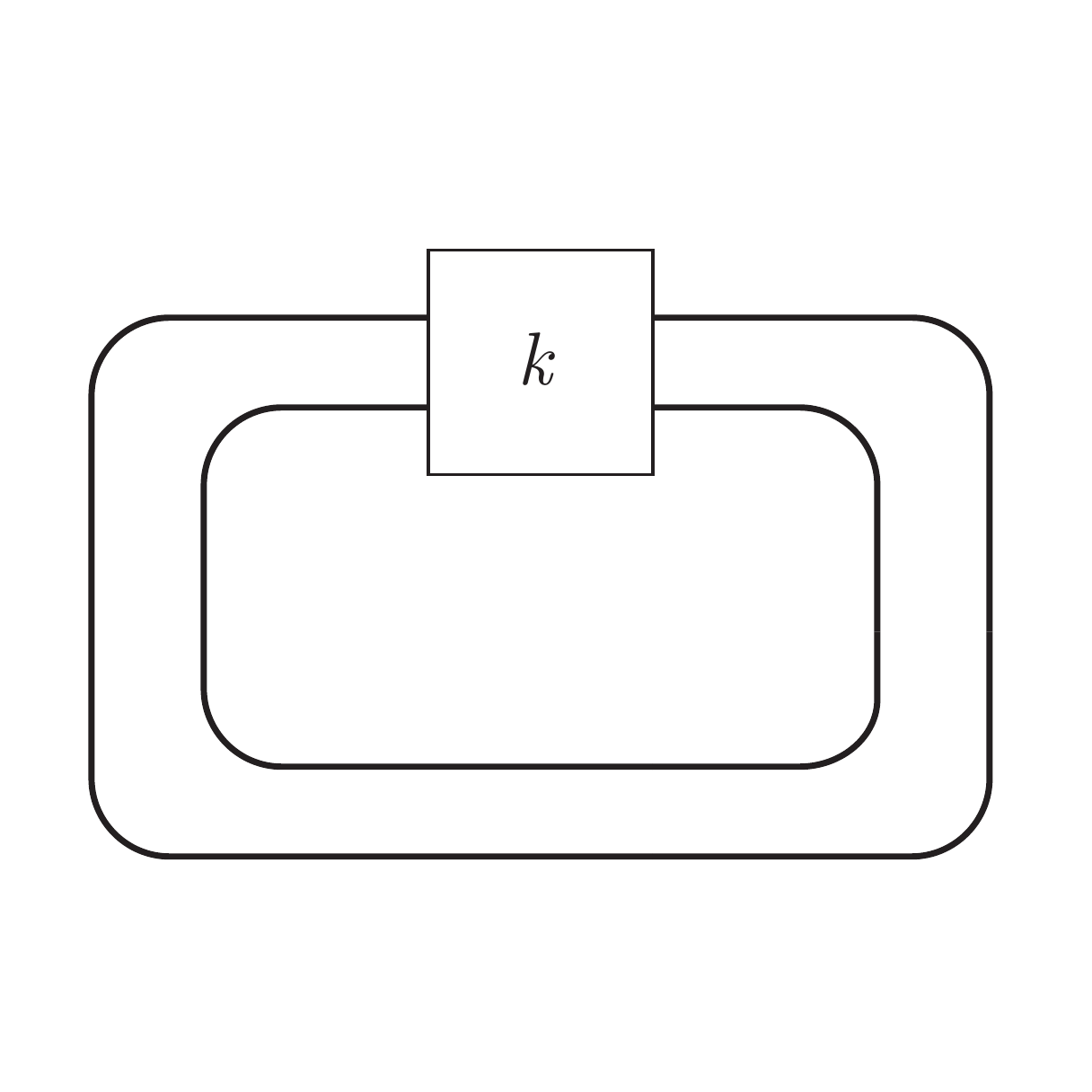} \\
\scriptstyle n =3 & &\scriptstyle n=2& &\scriptstyle k\neq0,\ n=2 \\
\\
\includegraphics[width=1.75in,trim=0 100pt 0 65pt,clip]{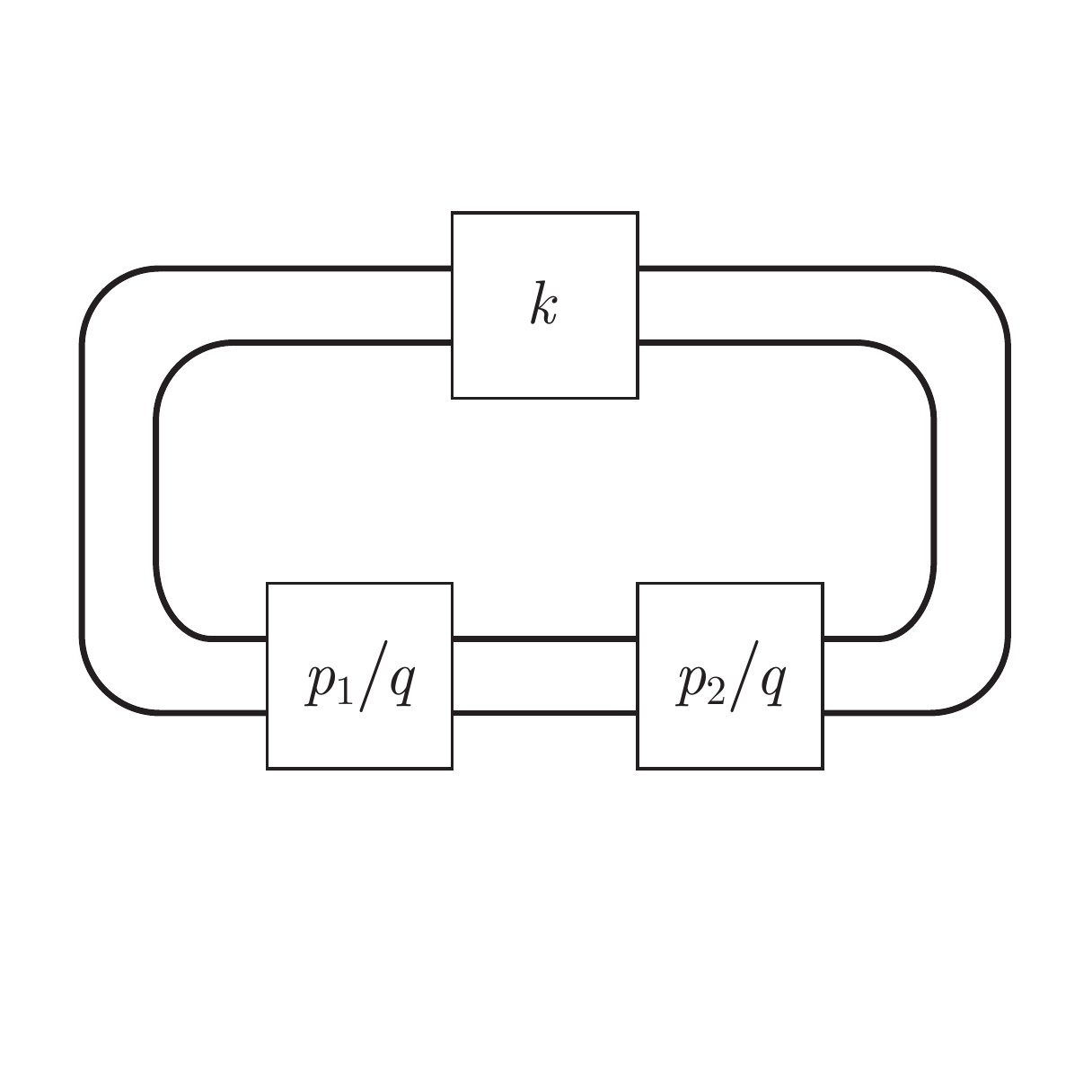} & & \includegraphics[width=1.65in,trim=0 85pt 0 65,clip]{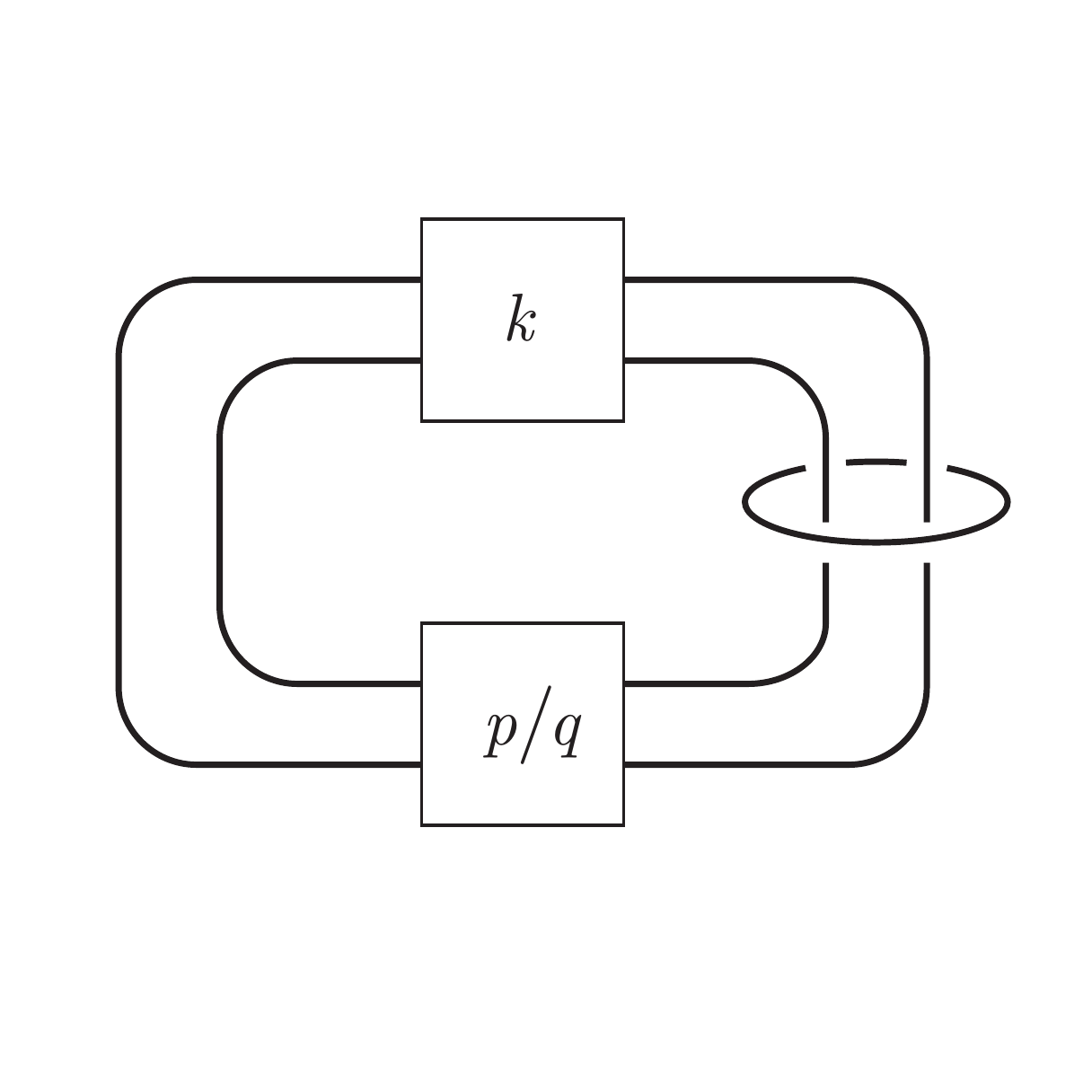} & & \includegraphics[width=2.00in,trim=0 100pt 0 75pt,clip]{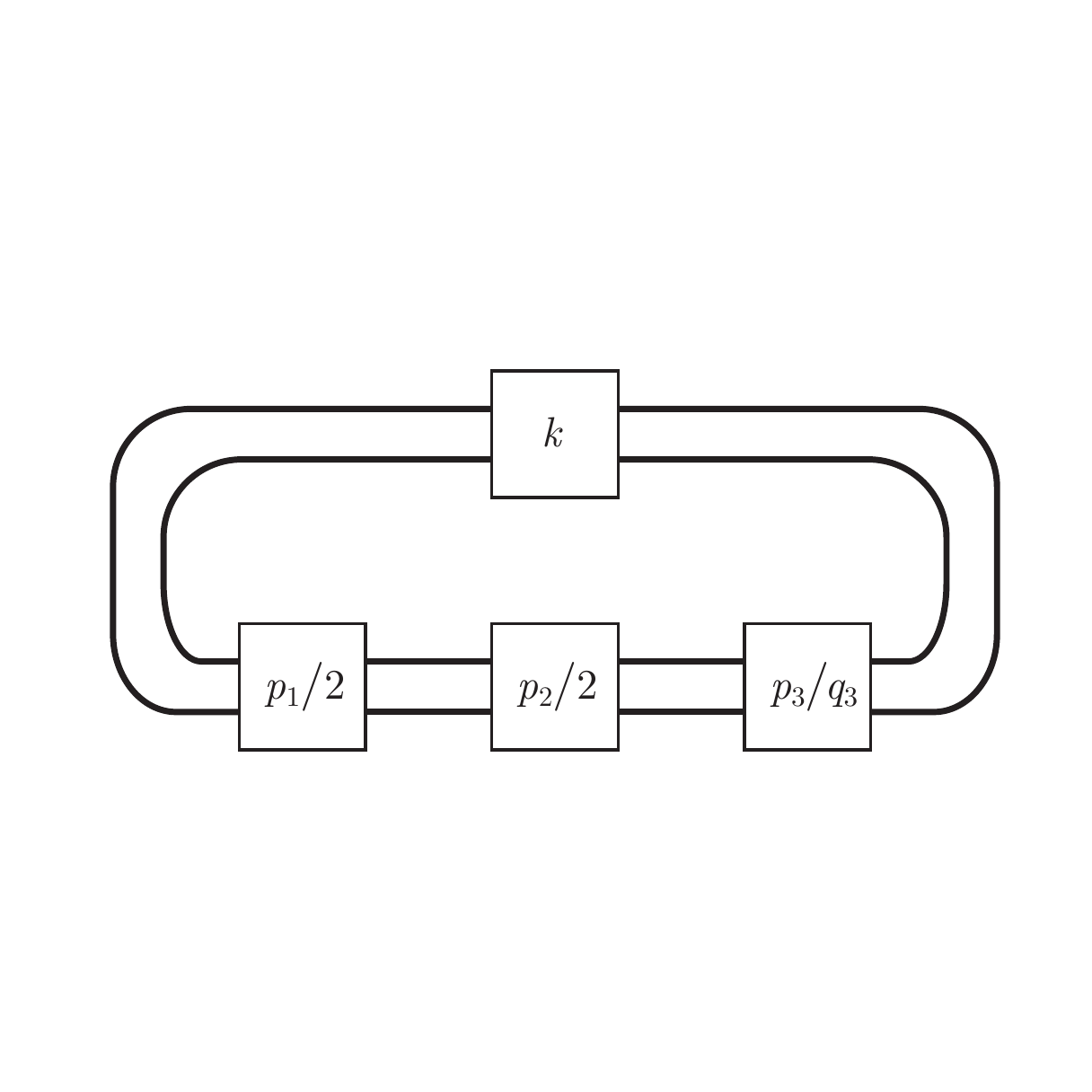} \\
\scriptstyle k+p_1/q+p_2/q \neq 0,\ n =2 & &\scriptstyle n=2& &\scriptstyle k+p_1/2+p_2/2+p_3/q_3 \neq 0,\ n =2 \\
\\
\includegraphics[width=2.00in,trim=0 90pt 0 75pt,clip]{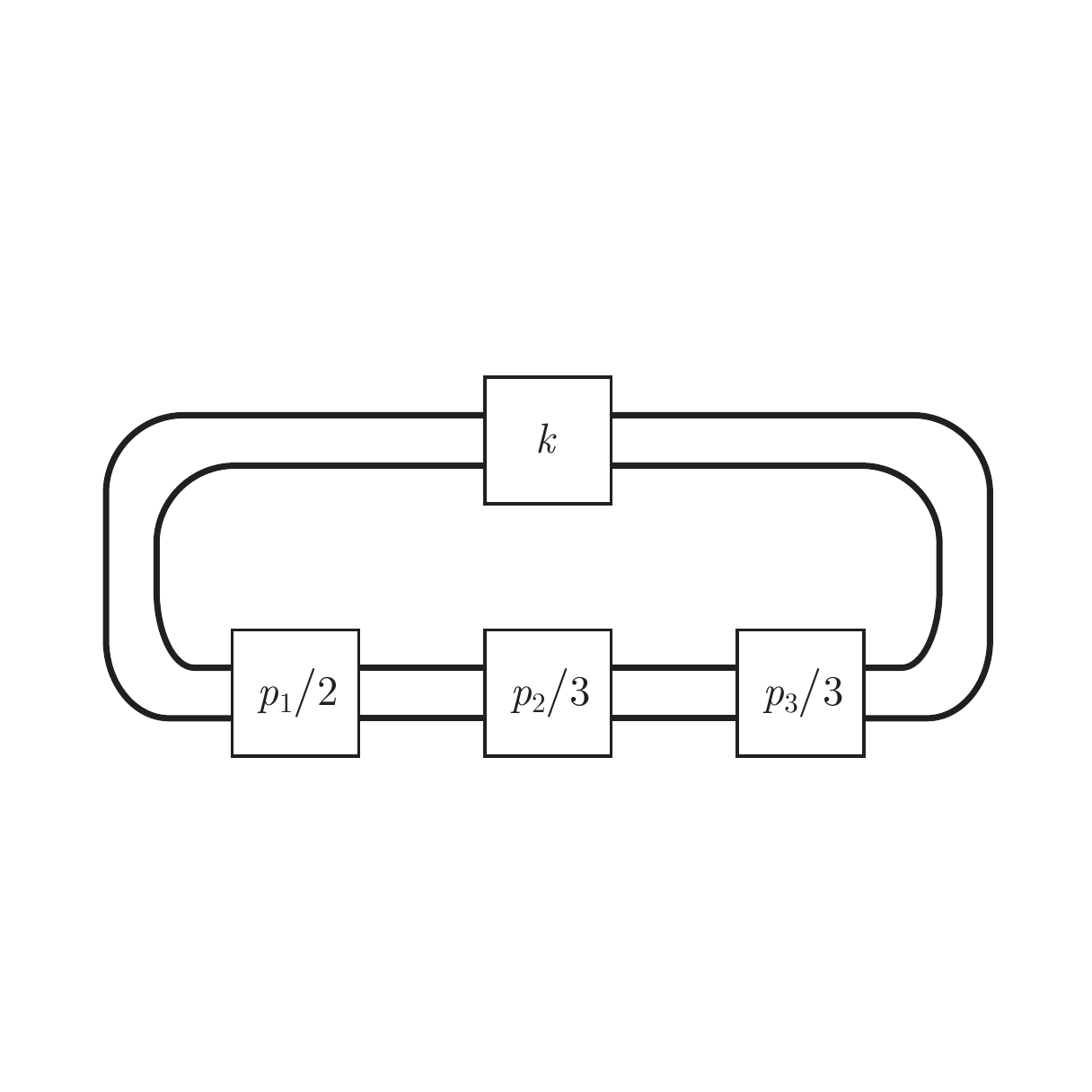} & & \includegraphics[width=2.00in,trim=0 100pt 0 75pt,clip]{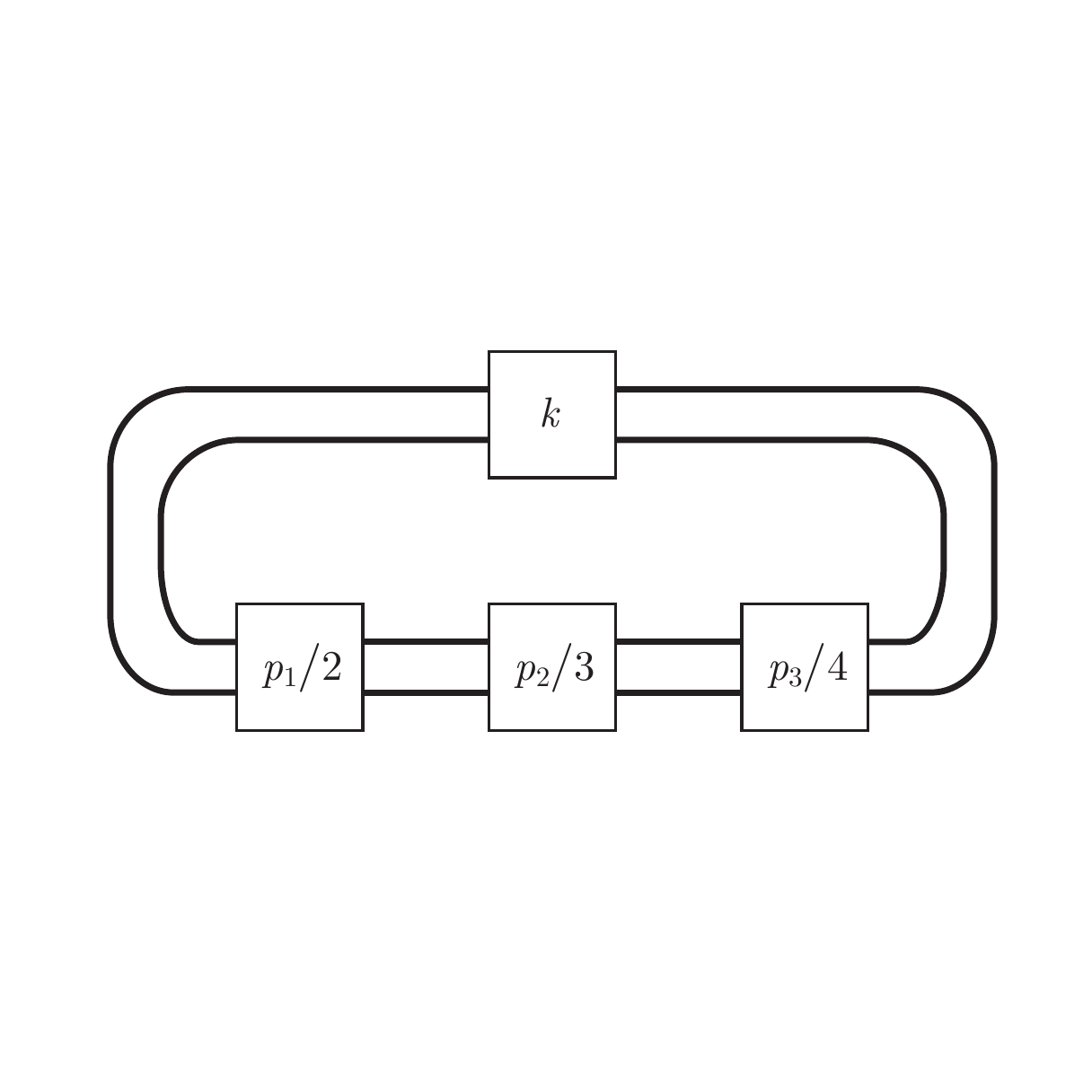} & & \includegraphics[width=2.00in,trim=0 100pt 0 75pt,clip]{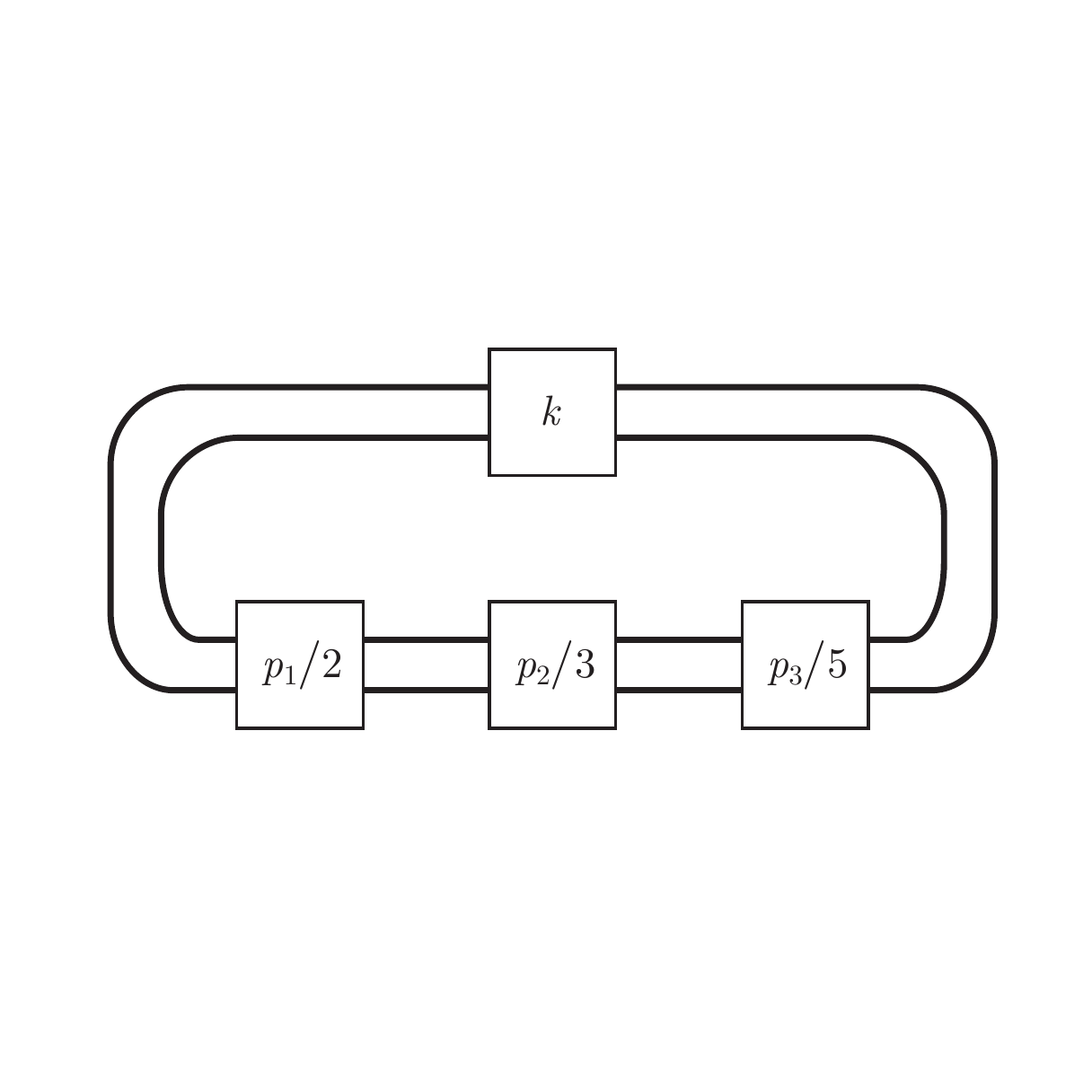}\\
\scriptstyle k+p_1/2+p_2/3+p_3/3 \neq 0,\ n =2 & & \scriptstyle k+p_1/2+p_2/3+p_3/4 \neq 0,\ n =2 & & \scriptstyle k+p_1/2+p_2/3+p_3/5 \neq 0,\ n =2
\end{array}
$$
\caption{Links $L \subset \mathbb S^3$ with finite $Q_n(L)$}
\end{table}

 \end{document}